\documentclass[aop,MSNbibl,seceqn,dvips]{arximspdf}

\doi{10.1214/12-AOP827} 
\volume{42}
\issue{3}
\pubyear{2014}
\firstpage{1121}
\lastpage{1160}

\makeatletter
\newproclaim{remark}{Remark}
\newproclaim{rems}{Remarks}
\newcommand{\rrvert}{\vert}
\newcommand{\llvert}{\vert}
\newtheorem{theorem}{Theorem}
\newtheorem{theo}{Theorem}[section]
\newproclaim{defi}{Definition}[section]
\newtheorem{propo}{Proposition}[section]
\newtheorem{lem}{Lemma}[section]

\newcommand{\PP}{\mathbb{P}}
\newcommand{\R}{\mathbb{R}}
\newcommand{\E}{\mathbb{E}}
\makeatother

\begin{document}
\begin{frontmatter}

\title{Cluster sets for partial sums and partial sum~processes}
\runtitle{Cluster sets}

\begin{aug}
\author[a]{\fnms{Uwe} \snm{Einmahl}\corref{}\thanksref{t1}\ead[label=e1]{ueinmahl@vub.ac.be}}
\and
\author[b]{\fnms{Jim} \snm{Kuelbs}\ead[label=e2]{kuelbs@math.wisc.edu}}
\thankstext{t1}{Supported in part by FWO-Vlaanderen Grant 1.5.167.09.}
\runauthor{U. Einmahl and J. Kuelbs}
\affiliation{Vrije Universiteit Brussel and University of Wisconsin}
\address[a]{Department of Mathematics\\
Vrije Universiteit Brussel\\
Pleinlaan 2\\
B-1050 Brussel\\ Belgium\\
\printead{e1}}

\address[b]{Department of Mathematics\\
University of Wisconsin\\
Madison, Wisconsin 53702\\
USA\\
\printead{e2}}
\end{aug}

\received{\smonth{4} \syear{2012}}
\revised{\smonth{12} \syear{2012}}

%
\begin{abstract}
Let $X, X_1, X_2,\ldots$ be i.i.d. mean zero random vectors with values
in a separable Banach space $B$, $S_n=X_1+\cdots+X_n$ for $n \ge1$,
and assume $\{c_n\dvtx  n\ge1\}$ is a suitably regular sequence of
constants. Furthermore, let
$S_{(n)}(t), 0 \le t \le1$ be the corresponding linearly interpolated
partial sum processes. We study the cluster sets $A= C(\{S_n/c_n\})$
and $\mathcal{A}=C(\{S_{(n)}(\cdot)/c_n\})$. In particular, $A$ and
$\mathcal{A}$ are shown to be nonrandom, and we derive criteria when
elements in $B$ and continuous functions $f\dvtx  [0,1] \to B$ belong to $A$
and $\mathcal{A}$, respectively. When $B= \mathbb{R}^d$ we refine our
clustering criteria to show both $A$ and $\mathcal{A}$ are compact,
symmetric, and star-like, and also obtain both upper and lower bound
sets for $\mathcal{A}$. When the coordinates of $X$ in $\mathbb{R}^d$
are independent random variables, we are able to represent $\mathcal
{A}$ in terms of $A$ and the classical Strassen set $\mathcal{K}$, and,
except for degenerate cases, show $\mathcal{A}$ is strictly larger than
the lower bound set whenever $d \ge2$. In addition, we show that for
any compact, symmetric, star-like subset $A$ of $\mathbb{R}^d$, there
exists an
$X
$ such that the corresponding functional cluster set $\mathcal{A}$ is
always the lower bound
subset. If $d=2$, then
additional refinements identify $\mathcal{A}$ as a subset of $\{
(x_1g_1, x_2g_2)\dvtx  (x_1,x_2) \in A, g_1,g_2 \in\mathcal{K}\}$, which is
the functional cluster set obtained when the coordinates are assumed to
be independent.
\end{abstract}

%


\end{frontmatter}
%
\section{Introduction}\label{sec1}
Let $X, X_1, X_2, \ldots$ be i.i.d. $d$-dimensional random vectors, and
let $S_n:=\sum_{j=1}^n X_j, n \ge1$. Denote the Euclidean norm on $\R
^d$ by $|\cdot|$ and write $\operatorname{cl}(M)$ for the closure of a subset
$M$ of a topological space.

Assuming $\E|X|^2 < \infty$ and
$\E X=0$, it follows from the $d$-dimensional version of the
Hartman--Wintner LIL that with
probability one,
%
%
\begin{equation}
\label{LIL} \limsup_{n \to\infty} |S_n|/\sqrt{2n\log\log n}
=\sigma,
\end{equation}
where $\sigma^2$ is the largest eigenvalue of the covariance matrix of
$X$.

For a sequence $\{x_n\dvtx n \ge1\} \subset\R^d$ the set of its limit
points is given by $\bigcap_{m=1}^{\infty} \operatorname{cl}(\{x_n\dvtx n \ge
m\}
)$. We denote this set by $C(\{x_n\dvtx n \ge1\})$, and we call it the
\textit{cluster set} of this sequence. This is obviously a closed subset
of $\R^d$.

Equation~(\ref{LIL}) then implies that with probability one, $C(\{
S_n/\sqrt{2n\log\log n}\dvtx\break  n \ge3\})$ is a compact subset of the
Euclidean ball with center 0 and radius $\sigma$ which must contain at
least one point from the boundary of this ball.

It is known for sums of i.i.d. random vectors and for any sequence $c_n
\nearrow\infty$ that the cluster set $C(\{S_n/c_n\dvtx n \ge1\})$ is
deterministic; see \cite{Kue}. So if $c_n$ is a sequence of that type
such that with probability one,
\[
\limsup_{n \to\infty}|S_n|/c_n < \infty,
\]
we have $C(\{S_n/c_n\dvtx n \ge1\}) = A$ with probability one, where $A$
is a nonempty compact subset of $\R^d$.

It is an interesting question to determine the cluster sets in such
cases. In the classical setting considered above it is well known that
$A =\{\Sigma x\dvtx |x| \le1\}$, where $\Sigma$ is the unique positive
semi-definite symmetric matrix satisfying $\Sigma^2$ = covariance
matrix of $X$.

A number of authors have investigated when one has LIL-type results for
random vectors $X$ with $\E|X|^2=\infty$. We mention the work of Kuelbs
\cite{Kue2} which implies among other things that if $X$ is a mean zero
random vector such that $S_n/a_n$ converges in distribution to a
$d$-dimensional normal distribution, one has for the normalizing
sequence $c_n = a_{[2n/\log\log n]}\log\log n, n \ge3$, and for
$\sigma^2$ equal to the largest eigenvalue of the covariance matrix of
the limit distribution, with probability one,
\[
\limsup_{n \to\infty} |S_n|/c_n = \sigma
\]
if and only if
$\sum_{n=1}^{\infty} \PP\{|X| \ge c_n\} < \infty$.

Moreover, the cluster set $C(\{S_n/c_n\dvtx n \ge1\})$ is in this case
again equal to $\{\Sigma x\dvtx |x| \le1\}$, with $\Sigma$ being chosen so
that $\Sigma^2$ is equal to the covariance matrix of the limit
distribution of $S_n/a_n$. It is easy to see that this result implies
the $d$-dimensional Hartman--Wintner LIL (just choose $a_n =\sqrt{n}$)
so that this is an extension of (\ref{LIL}).

This last result was generalized in \cite{E-3} where an
infinite-dimensional version of the Klass LIL \cite{Kla} is given. The
normalizing sequence $\gamma_n$ used in this result specializes in the
domain of attraction case to $\sigma a_{[2n/\log\log n]}\log\log n$,
but can also be applied for certain random vectors which are not in the
domain of attraction of a normal distribution. In these cases it was
not clear at all what the cluster sets $C(\{S_n/\gamma_n\dvtx n \ge1\})$
could be, given that there is no limit distribution with covariance
matrix available.

In \cite{E-1} it was shown that the cluster sets for this result have
to be subsets of the Euclidean unit ball which are star-like and
symmetric with respect to 0. Somewhat surprisingly, it also turned out
that any closed set of this type which contains a vector $a$ with
$|a|=1$ actually occurs as a cluster set.

Furthermore, it was shown in \cite{E-1} that if $X=(X^{(1)},\ldots,
X^{(d)}) $ and the variables $X^{(1)},\ldots, X^{(d)}$ are independent,
then the cluster sets are from the subclass of sets which are the
closures of at most countable unions of standard ellipsoids. Moreover
all sets of this type also occur as cluster sets in this case. Here we
call an ellipsoid ``standard'' if the main axes coincide with the
coordinate axes. Another way to say this is that a standard ellipsoid
is a set of the form $\{D x\dvtx |x| \le1\}$ where $D$ is a diagonal
matrix.

The following result follows from Theorem 4.1 in Einmahl and Li \cite
{EL} noticing that condition
(\ref{cn3}) below for $\mathbb{R}^d$ valued random vectors implies
that $\beta_0$ is equal to $0$ in this theorem.

\renewcommand{\thetheorem}{\Alph{theorem}}
\begin{theorem}\label{thmA}
Let $X, X_1, X_2, \ldots$ be i.i.d. mean zero
random vectors, and let $\{c_n\}$ be sequence of positive constants
such that
%
%
\begin{equation}
\label{cn1} c_n/\sqrt n\nearrow\infty,
\end{equation}
and for every $\varepsilon>0$ there exists an $m_{\varepsilon} \geq1$ such that
%
%
\begin{equation}
\label{cn2} c_n/c_m \leq(1+\varepsilon)n/m\qquad \mbox{for }
m_{\varepsilon} \leq m <n.
\end{equation}
Assume further that
%
%
\begin{equation}
\label{cn3} \sum_{n=1}^{\infty} \PP\bigl\{|X| \ge
c_n\bigr\} < \infty.
\end{equation}
Then we have for
the sums $S_n=\sum_{j=1}^n X_j, n \ge1$ with probability one,
%
%
\begin{equation}
\limsup_{n \to\infty} |S_n|/c_n =
\alpha_0,
\end{equation}
where
%
%
\begin{equation}
\label{cn4} \alpha_0 = \sup \Biggl\{\alpha\ge0\dvtx \sum
_{n=1}^{\infty} n^{-1}\exp \biggl(-
\frac{\alpha^2 c^2_n}{2nH(c_n)} \biggr) = \infty \Biggr\},
\end{equation}
with $H(t):=\sup\{ \mathbb{E}[\langle v,X \rangle^2I\{|X| \le t\}]\dvtx |v| \le1\}, t \ge0$.
\end{theorem}

All aforementioned LIL results and also the law of a very slowly
varying function (see Theorem 2 in \cite{EL-1}) follow from this
theorem.

The purpose of the present paper is to investigate whether there are
also general functional LIL-type results available in this case and
what the corresponding cluster sets are. In the 1-dimensional case this
question has been completely settled in~\cite{E} where it has been
shown that whenever $\alpha_0 < \infty$ and assumption (\ref{cn3}) is
satisfied, the functional LIL holds with cluster set $\alpha_0
\mathcal
{K}$, where $\mathcal{K}$ is the cluster set as in the Strassen LIL.

Much less is known in the multidimensional setting. We refer to \cite
{Kue2} where an infinite-dimensional functional LIL is established for
Banach space valued random vectors in the domain of attraction of a
Gaussian law. Nothing seems to be known---even in the
finite-dimensional case---for random vectors outside the domain of
attraction of a Gaussian law.
Given the complexity of the cluster sets $C(\{S_n/c_n\dvtx n \ge1\})$ in
this case, one cannot expect a simple answer as in the 1-dimensional setting.

\section{Statement of main results}\label{sec2} To formulate our results we need
somewhat more notation. Throughout, $X,X_1,X_2,\ldots$ are i.i.d. mean
zero random vectors, and except for the results of Section~\ref{sec3} they are
$\mathbb{R}^d$ valued. Let $C_d[0,1]$ be the continuous functions from
$[0,1]$ to $\R^d$ with sup-norm $\|f\|=\sup_{0 \le t \le1}|f(t)|, f
\in C_d[0,1]$.

The partial sum process $S_{(n)}\dvtx \Omega\to C_d[0,1]$
of order $n$ is defined by
\[
S_{(n)}(t)=S_{[nt]} + \bigl(nt - [nt]\bigr)X_{[nt]+1},\qquad 0
\le t \le1.
\]
The cluster set $C(\{S_{(n)}/c_n\dvtx n \ge1\})$ is defined as for sums,
that is, as the set of all limit points of the sequence $S_{(n)}/c_n$
in $C_d[0,1]$. We shall show (see Propositon \ref{pr3.1} below) that this set
is also deterministic.

Furthermore, we say the partial sum process
sequence $\{S_{(n)}(\cdot)\}$ converges and clusters compactly with
respect to a sequence
$c_n \nearrow\infty$ if we have that $C(\{S_{(n)}/c_n\dvtx n \ge1\})
=:\mathcal{A}$ is a compact subset of $C_d[0,1]$ and with probability one
$\lim_{ n \rightarrow\infty} d(S_{(n)}/c_n,\mathcal{A}) =0$,
where the distance between a function $f \in C_d[0,1]$ and $\mathcal
{A}$ is defined as
$d(f,\mathcal{A})= \inf_{g \in\mathcal{A}} \|f-g\|$. We write in this
case $\{S_{(n)}/c_n\} \leadsto\mathcal{A}$.

If $AC_{0}[0,1]$ denotes the absolutely continuous real valued
functions on $[0,1]$ which are zero when $t=0$, then for $ g \in C
[0,1]$, we define
\[
I(g)= %
\cases{ \displaystyle\int_0^1
\bigl(g'(s)\bigr)^2\,ds, &\quad  $\displaystyle g \in AC_{0}[0,1],
\int_0^1 \bigl(g'(s)
\bigr)^2\,ds<\infty$,
\vspace*{2pt}\cr
+\infty,&\quad$\mbox{otherwise.}$}
\]
One important fact about the $I$-functional is that it has a unique
minimum over closed balls. More precisely, suppose $g \in C[0,1]$ and
$\varepsilon>0$.
Then there exists a unique function, which we denote by $g_{\varepsilon}$,
such that $\Vert g-g_{\varepsilon}\Vert  \leq\varepsilon$ and
\[
I(g_{\varepsilon}) = \inf_{h: \Vert g-h\Vert  \leq\varepsilon} I(h).
\]
The existence of this minimum is well known, and details, as well as
further references, can be found in \cite{gr} and \cite{KLL}. Letting $
\mathcal{K }$ be the subclass of all functions in $AC_{0}[0,1]$ where
$I(g) \le1$, we get the cluster set in Strassen's functional LIL for
real-valued random variables.

If $\E|X|^2 < \infty$, the $d$-dimensional version of Strassen's
functional LIL applies which says that then with probability one,
%
%
\begin{equation}
\{S_{(n)}/\sqrt{2n\log\log n}\}\qquad \mbox{is relatively compact in }
C_d[0,1]
\end{equation}
and
%
%
\begin{equation}
\label{fclu} \mathcal{A}=C\bigl(\{S_{(n)}/\sqrt{2n\log\log n}\}\bigr) =
\Biggl\{\Sigma (f_1,\ldots, f_d)^t\dvtx \sum
_{i=1}^d I(f_i) \le1\Biggr\},
\end{equation}
where again $\Sigma$ is the positive semi-definite symmetric matrix
satisfying $\Sigma^2 =$ covariance matrix of $X$.

It is known that one can obtain the cluster sets $A=C(\{S_n/\sqrt {2n\log\log n}\})$ from~(\ref{fclu}) since
$A=\{f(1)\dvtx f \in\mathcal{A}\}$. Interestingly this implication can be
reversed. A small calculation shows that if the covariance matrix is
diagonal, we also have $\mathcal{A}= \{ x_1  \mathcal{K} \times
\cdots
\times x_d  \mathcal{K}\dvtx x =(x_1,\ldots,x_d) \in A\}$. This can also
be proved in general after replacing the canonical basis in $\R^d$ by
an orthonormal basis which diagonalizes the covariance matrix of $X$.

One might wonder whether a related phenomenon can be true if
$E|X|^2=\infty$. A~necessary condition for having $\mathcal{A}$ as in
the diagonal covariance matrix case would be that $A$ has an extended
symmetry property, namely $x =(x_1,\ldots,x_d) \in A \Rightarrow(\pm
x_1, \ldots, \pm x_d) \in A$
as one can choose functions $f_i \in\mathcal{K}$ with $f_i(1) =\pm1,
1 \le i \le d$.

So one might hope that the above result holds in general if $A$ has
this property. But it will turn out that this is not the case. For any
possible cluster set $A=C(\{S_n/c_n\dvtx n \ge1\})$, there exists a
distribution such that the functional cluster set is equal to the
smaller set $\{xg\dvtx  x \in A, g \in\mathcal{K}\}$ which only for very
special cases matches the function set above. This also shows that
relation (\ref{25}) in the subsequent Theorem \ref{theo21} gives an
optimal result.
%
\begin{theo} \label{theo21}
Let $X, X_1, X_2, \ldots$ be i.i.d. mean zero random vectors in $\R
^d$, and assume that
$
\sum_{n=1}^{\infty} \PP(|X|\geq c_n)<\infty,
$
where $\{c_n\}$ satisfies (\ref{cn1}) and (\ref{cn2}).
If $\alpha_0 =\limsup_{n \to\infty} |S_n|/c_n < \infty$,
we have with probability one,
%
%
\begin{equation}
\label{relcom} \{S_{(n)}/c_n\}\qquad  \mbox{is relatively compact in $C_d[0,1]$}.
\end{equation}
Consequently, the cluster set $\mathcal{A}= C(\{S_{(n)}(\cdot)/c_n\dvtx n
\ge1\})$ is compact in $C_d[0,1]$.
Furthermore, we have
%
%
\begin{equation}
\label{24} \mathcal{A} \subset\alpha_{1}\mathcal{K} \times\cdots\times
\alpha _{d}\mathcal{K},
\end{equation}
where $ \alpha_{i} = \limsup_{n \to\infty} |S_n^{(i)}|/c_n, 1 \le i
\le d$
and
%
%
\begin{equation}
\label{25} \mathcal{A} \supset\{ x g\dvtx x \in A, g \in\mathcal{K}\},
\end{equation}
where $A=C(\{S_n/c_n\dvtx n \ge1\}) \subset\R^d$.

Finally, $\mathcal{ A}$ is star-like and symmetric with respect to zero.
If
$f \in\mathcal{ A}$, then $f\dvtx [0,1] \rightarrow A$ continuously and for
$0 \leq t \leq1,  f(t) \in\sqrt tA$.
\end{theo}
%
\begin{remark*}
Using once more the fact that
$A=\{f(1)\dvtx f \in\mathcal{A}\}$,
we can conclude that these cluster sets are compact subsets of $\R^d$,
which are star-like and symmetric with respect to zero. This has been
proven in \cite{E-1} only for a special case of Theorem~\ref{thmA}. We now see
that this is always the case when the assumptions of Theorem~\ref{thmA} are
satisfied.
\end{remark*}

If the coordinates $X^{(1)},\ldots, X^{(d)}$ of $X$ are independent,
our next result
gives the complete answer showing that in this case we again have a
1--1 correspondence between the functional cluster sets $\mathcal
{A}=C(\{
S_{(n)}/c_n\dvtx n \ge1\})$ and $A=C(\{S_n/c_n\dvtx n \ge1\})$.
%
\begin{theo} \label{theo22}
Let $X =(X^{(1)},\ldots,X^{(d)})\dvtx \Omega\to\R^d$ be a mean zero random
vector with independent components
and suppose that $\{c_n\}$ satisfies (\ref{cn1}) and (\ref{cn2}).
If $\alpha_{0} < \infty$ and $\sum_{n=1}^{\infty} \PP(|X|\geq
c_n)<\infty$, then
with probability one we have (\ref{relcom}) and
$\mathcal{A}= \{ x_1  \mathcal{K} \times\cdots\times x_d
\mathcal
{K}\dvtx x =(x_1,\ldots,x_d) \in A\}$.
\end{theo}
Interestingly it turns out that if $d=2$, the last set is also the
maximal set for the cluster sets in the general case. It is not clear
whether this is also the case in higher dimensions.
%
\begin{theo} \label{theo23}
Let $X$ be a mean zero random vector in $\R^2$, and assume that
$
\sum_{n=1}^{\infty} \PP(|X|\geq c_n)<\infty,
$
where $\{c_n\}$ satisfies (\ref{cn1}) and (\ref{cn2}).
If $\alpha_0 < \infty$, we always have $\mathcal{A} \subset\{ x_1
\mathcal{K} \times x_2  \mathcal{K}\dvtx x \in A\}$.
\end{theo}
The remaining part of the paper is organized as follows: In Section~\ref{sec3}
we prove some general results on cluster sets in the functional LIL.
Though the present paper considers mainly the finite-dimensional case
we establish these results in the infinite-dimensional setting so that
they can be used in future work on the functional LIL problem in this
more general setting. In Section~\ref{sec4} we then derive via a strong
approximation result of Sakhanenko \cite{sakh} criteria for clustering in
$\R
^d$ in terms of Brownian motion probabilities. This enables us in
Sections~\ref{sec5}--\ref{sec7} to prove Theorems \ref{theo21}, \ref{theo22} and \ref
{theo23} using results on Gaussian probabilities of balls in
$(C_d[0,1], \|\cdot\|)$.
Finally, in Section~\ref{sec8} we shall provide an example where the cluster set
$A=C(\{S_n/c_n\dvtx n \ge1\})$ is equal to an arbitrary given closed,
star-like and symmetric set $\tilde{A}$ with $\max_{x \in\tilde{A}}
|x|=1$ and at the same time the functional cluster set $\mathcal{A}$ is
equal to $\{xg\dvtx  x \in\tilde{A}, g \in\mathcal{K}\}$.

\section{Some general results on cluster sets}\label{sec3}
Here we present results for the cluster sets $C(\{S_{(n)}/c_n\dvtx n \ge1\}
)$ and $C(\{S_n/c_n\dvtx n \ge1\})$. They include their behavior when the
sequences $\{S_{(n)}/c_n\}$ and $\{S_n/c_n\}$ are relatively compact
with probability one. Moreover, we provide a necessary and sufficient
series condition characterizing the functions $f$ in the functional
cluster sets $C(\{S_{(n)}/c_n\dvtx n \ge1\})$. As our proofs work also in
the infinite-dimensional setting, we now consider $B$-valued random
variables $X, X_1, X_2, \ldots,$ where $(B,|\cdot|)$ is a separable
Banach space with norm $|\cdot|$.
\subsection{Nonrandomness of the functional cluster sets}\label{sec3.1}
Our first result is a zero-one law showing the cluster set $C(\{
S_{(n)}/c_n\dvtx n \ge1\})$ is deterministic with probability one, and is
the analogue of Lemma 1 in \cite{Kue}.

Let for $0 \le m \le n$
\begin{eqnarray}
S_{(n,m)}(t)= %
\cases{ 0, &\quad$0 \le t \le m/n$,
\vspace*{2pt}\cr
S_k - S_m, &\quad $t=k/n, m \le k \le n$,
\vspace*{2pt}\cr
\mbox{linearly interpolated elsewhere},}\nonumber\\
\eqntext{(0 \le t \le1).}
\end{eqnarray}
Obviously, the choice $m=0$ gives us the partial sum process $S_{(n)}$
of order $n$, and these processes are random elements in the space
$C_0([0,1],B)$ of all continuous functions $f\dvtx [0,1] \to B$ satisfying
$f(0)=0$. We denote the sup-norm on this space by $\|\cdot\|$.

\begin{propo}\label{pr3.1} Let $\{c_n\}$ be a positive sequence such that $c_n
\nearrow\infty$. Then, there exists a nonrandom set $\mathcal{A}$ in
$C_0([0,1],B)$ such that with probability one
%
%
\begin{equation}
\label{cluster1} C\bigl(\bigl\{S_{(n)}(\cdot)/c_n\dvtx n
\ge1\bigr\}\bigr) = \mathcal{A}.
\end{equation}
\end{propo}

\begin{pf} First of all observe that the Banach space $C_0([0,1],B)$
is separable. This follows since $B$ separable implies one can embed
$B$ into $C[0,1]$. Then $C_0([0,1],B)$ is embedded isometrically into
$C_0([0,1],C[0,1])$, and the polynomials in two variables and rational
coefficients are dense in this space. Hence there exists a countable
base $\mathcal{B}$ for the norm topology of $C_0([0,1],B)$.

Let
\[
\mathcal{B}_1 = \Bigl\{U \in\mathcal{B}\dvtx \PP\Bigl(\liminf
_{n \to\infty} d(S_{(n)}/c_n, U)=0\Bigr)=1\Bigr\}
\]
and
\[
\mathcal{B}_2 =\Bigl\{U \in\mathcal{B}\dvtx \PP\Bigl(\liminf
_{n \to\infty} d(S_{(n)}/c_n, U)=0\Bigr)=0\Bigr\}.
\]
As we have for any fixed $m$ and $n \ge m$, $\|S_{(n)} - S_{(n,m)}\|
/c_n \to0$ as $n \to\infty$, we see that
\[
\Bigl\{\liminf_{n \to\infty} d(S_{(n)}/c_n, U)=0
\Bigr\}=\Bigl\{\liminf_{n \to
\infty} d(S_{(n, m)}/c_n,
U)=0\Bigr\},\qquad U \in\mathcal{B}.
\]
The event on the right-hand side is measurable with respect to the
$\sigma$-field generated by $X_{m+1}, X_{m+2}, \ldots$ and this holds
for any $m$.

Thus $\{\liminf_{n \to\infty} d(S_{(n)}/c_n, U)=0\}$ is
a tail event, and
by Kolmogorov's zero one law we have
$\mathcal{B} =\mathcal{B}_1 \cup\mathcal{B}_2$.\vadjust{\goodbreak}

Let $V=\bigcup_{U \in\mathcal{B}_2} U$ and $\mathcal{A}= C_0([0,1],B)
\setminus V$. Then $\mathcal{A}$ is nonrandom, and we set
\[
\Omega_1= \bigcap_{U \in\mathcal{B}_1} \Bigl\{\omega
\dvtx \liminf_{n \to
\infty
} d\bigl(S_{(n)}(\omega,
\cdot)/c_n, U\bigr)=0\Bigr\}
\]
and
\[
\Omega_2= \bigcap_{U \in\mathcal{B}_2}\Bigl\{\omega
\dvtx \liminf_{n \to
\infty} d\bigl(S_{(n)}(\omega,
\cdot)/c_n, U\bigr)>0\Bigr\}.
\]
Then $\Omega_1 \cap\Omega_2$ is the countable intersection of sets of
probability one, so it has probability one.
So it is sufficient to prove that we have for $\omega\in\Omega_1
\cap
\Omega_2$,
%
%
\begin{equation}
\label{cluster2} D(\omega)\equiv\bigcap_{m=1}^{\infty}
\operatorname{cl}\bigl(\bigl\{S_{(n)}(\omega,\cdot)/c_n\dvtx n
\geq m\bigr\}\bigr)=\mathcal{A}.
\end{equation}
To prove (\ref{cluster2}), we first note that for $g \in\mathcal{A}$
and $\varepsilon>0$ there is a $U \in\mathcal{B}$ with $g \in U \subset
U_{\varepsilon}(g)$, where as usual $U_{\varepsilon}(g) =\{f \in
C_0([0,1],B)\dvtx \|f- g\| < \varepsilon\}$. As $g \notin V$, this implies
$U \notin\mathcal{B}_2$ so $U \in\mathcal{B}_1$. Hence by
definition of $\Omega_1$
we have $S_{(n)}(\omega,\cdot)/c_n \in U_{2\varepsilon}(g)$ infinitely
often. Therefore, since $\varepsilon$ is arbitrary, $g \in D(\omega)$ and
hence $\mathcal{A} \subset D(\omega)$ for all $ \omega\in\Omega_1
\cap\Omega_2$.

On the other hand, if $g \notin\mathcal{A}$ or equivalently, $g \in
V$ there is a $U \in\mathcal{B}_2$ with $g \in U$. By definition of
$\Omega_1 \cap\Omega_2$ we have $S_{(n)}(\omega,\cdot)/c_n \in U^c$
eventually. Hence $ g \notin D(\omega)$, and therefore $D(\omega)
\subset\mathcal{A}$ and (\ref{cluster2}) has been proven.
\end{pf}
\subsection{Compactness of the functional cluster sets}\label{sec3.2}
%
\begin{propo}\label{pr3.2} Let $\{c_n\}$ be a positive sequence such that $c_n
\nearrow\infty$, and assume $\mathcal{A}$ is the deterministic cluster
set of $S_{(n)}/c_n$ determined as in (\ref{cluster1}).
If $\{S_{(n)}/c_n\}$ is relatively compact in $C_0([0,1],B)$ with
probability one, then $\mathcal{A}$ is a compact nonempty subset of
$C_0([0,1],B)$ and with probability one $S_{(n)}/c_n$ converges and
clusters compactly to $\mathcal{A}$, that is, with probability one
$\{S_{(n)}/c_n\} \leadsto\mathcal{A}$.
\end{propo}
\begin{pf} Let $\mathcal{A}$ be the deterministic cluster set of $\{
S_{(n)}/c_n\}$. We claim that $\{S_{(n)}/c_n\}$ relatively compact in
$C_0([0,1],B)$ with probability one implies
$\lim_{n \rightarrow\infty} d(S_{(n)}/c_n,\mathcal{A})=0$
with probability one.

To see this, suppose that
$
\limsup_{n \rightarrow\infty} d(S_{(n)}/ c_n,\mathcal{A})>0
$
with positive probability. Then there is a $\delta>0$ such that with
positive probability\break  $\limsup_{n \rightarrow\infty} d(S_{(n)}/ c_n,\mathcal{A}) \ge
2\delta.
$
Now the set $E=\{x\dvtx d(x,\mathcal{A}) \ge\delta\}$
is closed, and with positive probability the relatively compact
sequence $\{S_{(n)}/c_n\}$ would be infinitely often in $E$ and would
have limit points in $E$ which is impossible since $\mathcal{A} \cap E
= \varnothing$.

Finally, $\mathcal{A}$ is compact and nonempty as $\mathcal{A}=
\bigcap_{m\geq1}\operatorname{cl}(\{S_{(n)}( \omega)/c_n\dvtx\break  n \geq m\})$ with
probability one. Choosing $\omega$ so that this holds and at the same
time $\operatorname{cl}(\{S_{(n)}( \omega)/c_n\dvtx n \geq1\})$ is compact, we
readily obtain that the closed set $\mathcal{A}$ is compact as
well.\vadjust{\goodbreak}
\end{pf}

Our next proposition relates the clustering and compactness of $\{
S_n/c_n\}$ to the clustering and compactness of $\{S_{(n)}(\cdot)/c_n\}$
in Banach spaces where one has finite rank operators that approximate
the identity. More precisely, a Banach space $B$ has the approximation
property if for each compact subset $K$ of $B$ and $\varepsilon>0$ there
is a finite rank operator $T\dvtx B \rightarrow B$ such that
\[
\sup_{x \in K}\bigl|x -T(x)\bigr| < \varepsilon.
\]
This property is less restrictive than requiring $B$ have a Schauder
basis, and hence many (but not all) Banach spaces have the
approximation property. Information about this property is easily
found, and two classical references are \cite{Day} and \cite{LT}.
%
\begin{propo}\label{pro33}
Let $\{c_n\}$ satisfy (\ref{cn1}) and (\ref{cn2}), and assume
%
%
\begin{equation}
\label{f35} \sum_{n= 1}^{\infty} \PP\bigl(|X|
>c_n\bigr) < \infty.
\end{equation}
If $(B,|\cdot|)$ has the approximation property and $\{S_n/c_n\}$ is
relatively compact in~$B$ with probability one, then $\{S_{(n)}(\cdot
)/c_n\}$ is relatively compact in $C_0([0,1],B)$ with probability one.
Moreover, if $\mathcal{A}$ is the deterministic cluster set for $\{
S_{(n)}(\cdot)/c_n\}$ given in (\ref{cluster1}), then $\mathcal{A}$ is
nonempty and compact and we have with probability one, $\{
S_{(n)}/c_n\} \leadsto\mathcal{A}$.
\end{propo}
\begin{pf} To verify this let $\varepsilon>0$ be given. Since $\{
S_n/c_n\}
$ is relatively compact in $B$ with probability one, then by the same
argument as in Proposition \ref{pr3.2} the deterministic cluster set $A$ of $\{
S_n\}$ with respect to $\{c_n\}$ is compact and such that with
probability one
\[
\PP\Bigl(\limsup_{n \rightarrow\infty} d(S_n/c_n,A)=0
\Bigr)=1
\]
and
\[
\PP\bigl(C\bigl(\{S_n/c_n\dvtx n \geq1\}\bigr)=A
\bigr)=1.
\]
Since $(B,|\cdot|)$ has the approximation property, given $\varepsilon>0$
there exists a finite rank operator
\[
\Lambda(x) =\sum_{i=1}^d
f_i(x)x_i
\]
mapping $B$ into $B$, with $x_1,\ldots,x_d \in B$ and $f_1,\ldots,f_d
\in B_1^{*}$, the unit ball of $B^*$, such that
\[
\sup_{x \in A} \bigl|x - \Lambda(x)\bigr| < \varepsilon.
\]
Then, with probability one
\[
\limsup_{n \rightarrow\infty} \bigl|S_n/c_n -
\Lambda(S_n/c_n)\bigr| \leq \varepsilon,
\]
and hence we also have
\[
\limsup_{n \rightarrow\infty} \sup_{0\leq t \leq1}\bigl|S_{(n)}(t)/c_n-
\Lambda\bigl(S_{(n)}(t)/c_n\bigr)\bigr| \leq\varepsilon
\]
with probability one.

Now let
%
%
\begin{equation}
\sigma_{n,i}^2 = E\bigl(f_i(X)^2I
\bigl(\bigl|f_i(X)\bigr| \leq c_n\bigr)\bigr),\qquad  1 \leq i \leq d,
\end{equation}
and define
%
%
\begin{equation}
\alpha_i =\sup\biggl\{ \alpha\geq0\dvtx \sum
_{n \geq1} n^{-1}\exp\biggl\{-\frac
{\alpha
^2c_n^2}{2n\sigma_{n,i}^2}\biggr\} =
\infty\biggr\}
\end{equation}
for $i=1,\ldots,d$. Also let $\mathcal{K}$ denote the limit set in the
functional law of the iterated logarithm for Brownian motion as defined
in Section~\ref{sec2}.

Each random variable $f_i(X), i=1,\ldots,d$, is such that
$E(|f_i(S_n/c_n)|) \rightarrow0$ since the real line is a type 2 Banach
space. See Lemma 4.1 in \cite{EL}. In addition, since the $f_i$'s are continuous linear
functionals in $B_1^*$, and $S_n/c_n$ is relatively compact in $B$ with
probability one, we have
from (\ref{f35}) that for $i=1,\ldots,d$
\[
\sum_{n=1}^{\infty} \PP\bigl(\bigl|f_i(X)\bigr|
> c_n\bigr) < \infty,
\]
and with probability one
\[
\limsup_{n\rightarrow\infty}\bigl|f_i(S_n/c_n)\bigr|
< \infty,\qquad i=1, \ldots,d.
\]
Hence (4.4) of Theorem 5 of \cite{EL} implies with probability one that
\[
\limsup_{n \rightarrow\infty} \bigl|f_i(S_n/c_n)\bigr|
=\alpha_i,
\]
and since this $\limsup$ is finite with probability one we have $\alpha_i<
\infty, i=1,\ldots,d$.

Thus Theorem 1 of \cite{E} implies that for every $\varepsilon>0$
\[
\PP\Biggl(\bigcap_{i=1}^d\bigl
\{f_i\bigl(S_{(n)}(\cdot)/c_n\bigr) \in(
\alpha_i\mathcal {K})^{\varepsilon} \mbox { eventually}\bigr\}
\Biggr)=1,
\]
and hence by the equivalence of norms on finite dimensional Banach
spaces we also have
\[
\PP\bigl(\Lambda\bigl(S_{(n)}(\cdot)/c_n\bigr) \in(
\alpha_1\mathcal{K}\times \cdots \times\alpha_d
\mathcal{K})^{\varepsilon} \mbox { eventually}\bigr) =1
\]
for all $\varepsilon>0$.
Therefore, we have $\{S_{(n)}(\cdot)/c_n\dvtx n \geq1\}$ totally bounded,
and thus relatively compact, in $C_0([0,1],B)$ with probability one.
Proposition \ref{pr3.1} now implies $\mathcal{A}$ is a nonempty compact set
and that $\{S_{(n)}/c_n\} \leadsto\mathcal{A} $ with probability one.
\end{pf}
\subsection{\texorpdfstring{The functional LIL version of a result of Kesten \cite{Ke}}
{The functional LIL version of a result of Kesten [11]}}\label{sec3.3}
The purpose of this part of the paper is to derive a necessary and
sufficient condition that a function $f \in C_0([0,1],B)$ is in the
deterministic cluster set $\mathcal{A}=C(\{S_{(n)}/c_n\dvtx n \ge1\})$,
where we use the same notation as in Section~\ref{sec3.1}. The corresponding
result for the cluster set $A=C(\{S_n/c_n\dvtx n \ge1)$ (see~Lemma 1 in
\cite{E-1}) reads as follows:
%
%
\begin{eqnarray}
\label{Kes} \qquad x \in C\bigl(\{S_n/c_n\}\bigr) \quad\mbox{a.s.}\quad
\Longleftrightarrow\quad
\sum_{n=1}^{\infty}
n^{-1}\PP\bigl\{|S_n/c_n - x|< \varepsilon\bigr\} =
\infty,\quad \varepsilon > 0,
\nonumber
\\[-8pt]
\end{eqnarray}
where one has to assume that $S_n/c_n$ is stochastically bounded, and
$c_n$ satisfies conditions (\ref{cn1}) and (\ref{cn2}).
This result for real-valued random variables goes back to Theorem 3 in
Kesten \cite{Ke} who actually considers somewhat more general sequences
$\{c_n\}$.

We now prove such a result for partial sum processes based on i.i.d.
mean zero random variables taking values in a separable Banach space
$(B,|\cdot|)$. To simplify notation we set $s_n = S_{(n)}/c_n, n \ge
1$, and we denote the sup-norm of any continuous function $f\dvtx [0,1] \to
B$ by $\|f\|$.
%
\begin{propo} \label{pro1}
Let $f\dvtx [0,1] \to B$ continuous, and let $c_n$ be a sequence of
positive real numbers satisfying conditions (\ref{cn1}) and (\ref
{cn2}). Take a fixed $\rho> 1$. Then
the following are equivalent:
\begin{longlist}[(a)]
\item[(a)] $f \in C(\{s_n\dvtx n \ge1\})\ \mathrm{a.s.}$;
\item[(b)] $\sum_{k= 0}^{\infty} \PP\{\|s_n -f\| < \varepsilon\mbox{ for
some }n \in[\rho^k,\rho^{k+1}[ \} = \infty, \varepsilon> 0$.
\end{longlist}
\end{propo}
\begin{pf}
\fbox{(b) $\Rightarrow$ (a)} To further simplify our
notation, we set
$I_k = \{n\dvtx \rho^k \le n < \rho^{k+1}\}$ and
\[
G_k = \bigcup_{n \in I_k} \bigl\{\|s_n
- f\| < \varepsilon\bigr\}, \qquad k \ge0.
\]
Consider also the stopping times $\tau_k$ defined by
\[
\tau_k = \inf\bigl\{ n \ge\rho^k\dvtx \|s_n
- f\| < \varepsilon\bigr\}, \qquad k \ge0.
\]
Then we obviously have
%
%
\begin{equation}
G_k = \bigl\{\tau_k < \rho^{k+1}\bigr\},\qquad k
\ge0.
\end{equation}
Set
\[
H_k =\bigl\{\|s_n - f\| \ge\varepsilon\mbox{ for all }n
\ge\rho^{k+r}\bigr\} \cap G_k,
\]
where $ r > 0$ is an integer which will be specified later.

Then it is obvious that
%
%
\begin{equation}
\PP(H_k) = \sum_{m \in I_k} \PP\bigl\{
\|s_n - f\| \ge\varepsilon\mbox{ for all }n \ge\rho^{k+r},
\tau_k = m\bigr\}.
\end{equation}
Next set for $0 \le m \le n$,
$s_{n,m} = S_{(n,m)}/c_n$, where $S_{(n,m)}$ is defined as in
Section~\ref{sec3.1}.
Then we have for $m \in I_k$ and $n \ge\rho^{k+r}$ on the event $\{
\tau
_k = m\} \subset\{\|s_m -f\| < \varepsilon\}$
\begin{eqnarray*}
\|s_{n,m} - s_n\|& \le&\|S_{(m)}\|/c_n
\le\bigl (\|s_m -f\| + \|f\| \bigr)c_m/c_n
\\
&\le& \bigl(\varepsilon+ \|f\|\bigr)\sqrt{m/n}
\\
&\le& \bigl(\varepsilon+ \|f\|\bigr)\rho^{(1-r)/2} \le\varepsilon
\end{eqnarray*}
provided that we choose $r=r(\varepsilon, f)$ large enough.

Due to the independence of $s_{n,m}$ and the event $\{\tau_k =m\}$, we
can infer that
%
%
\begin{equation}
\PP(H_k) \ge\sum_{m \in I_k} \PP\bigl\{
\|s_{n,m} - f\| \ge2 \varepsilon \mbox { for all }n \ge
\rho^{k+r}\bigr\}\PP\{\tau_k = m\}.
\end{equation}
Next observe that
\[
S_{(n-m)}(t)_{0 \le t \le1} \stackrel{d} {=} S_{(n,m)}\bigl(
\alpha _{n,m}(t)\bigr)_{0 \le t \le1},
\]
where $\alpha_{n,m}(t) = (m/n)+(1-m/n)t, 0 \le t \le1$.

Set $f_{n,m}(t) = f(\alpha_{n,m}(t)), 0 \le t \le1$. Then it is easy
to see that by uniform continuity of $f$ we have $\|f - f_{n,m}\| <
\varepsilon$ if we have chosen $r$ large enough.
We conclude that
\begin{eqnarray*}
\PP(H_k) &\ge& \sum_{m \in I_k}\PP\bigl\{
\bigl\|c_n^{-1} S_{(n-m)} - f_{n,m}\bigr\| \ge2
\varepsilon\mbox{ for all } n \ge\rho^{k+r}\bigr\}\PP\{\tau_k
= m\}
\\
&\ge& \sum_{m \in I_k}\PP\bigl\{\bigl\|c_n^{-1}
S_{(n-m)} - f\bigr\| \ge3\varepsilon \mbox{ for all } n \ge\rho^{k+r}
\bigr\}\PP\{\tau_k = m\}.
\end{eqnarray*}
Moreover, for $\varepsilon>0$ and $f $ fixed, we take $\hat\varepsilon>0$
such that $\hat\varepsilon(\Vert f\Vert  \vee1)< \varepsilon$. Then, for large $k$
\[
\bigl\|f - (c_{n-m}/c_n)f\bigr\| \le\bigl(1 - (1+\hat
\varepsilon)^{-1}(1-m/n)\bigr)\|f\| = \Vert f\Vert (1+\hat
\varepsilon)^{-1}(\hat\varepsilon+m/n),
\]
which is $\le2\varepsilon$ if we choose $r$ large enough that $\Vert f\Vert \rho
^{1-r} < \varepsilon$.

Therefore,
\[
\PP(H_k) \ge\sum_{m \in I_k} \PP\bigl\{
\|s_{n-m} - f\| \ge5\varepsilon c_n/c_{n-m} \mbox{
for all }n \ge\rho^{k+r}\bigr\}\PP\{\tau_k =m\}.
\]
Assuming also that $r$ is so large that for sufficiently large $m$,
\[
c_n/c_{n-m} \le1.1 n/(n-m) \le1.2 \qquad\mbox{whenever } m/n
\le\rho^{1-r},
\]
we readily obtain from the last inequality
\[
\PP(H_k) \ge\sum_{m \in I_k} \PP\bigl\{
\|s_{n-m} -f\| \ge6\varepsilon \mbox{ for all }n \ge\rho^{k+r}
\bigr\} \PP\{\tau_k=m\},
\]
which in turn is
\[
\ge\PP\bigl\{\|s_n -f\| \ge6\varepsilon\mbox{ for all }n \ge
\rho^r-\rho \bigr\}\PP(G_k).
\]
Noticing that $H_k \cap H_{\ell} = \varnothing, |k-\ell| > r$, we see that
$ Y = \sum_{k=0}^{\infty} I_{H_k} \le r$
and consequently
%
%
\begin{equation}\qquad
r \ge\E[ Y] = \sum_{k=0}^{\infty}
\PP(H_k) \ge\PP\bigl\{\|s_n -f \| \ge6\varepsilon, n \ge
\rho^r-\rho\bigr\}\sum_{k=k_{\varepsilon}}^{\infty}
\PP(G_k).
\end{equation}
The last series is divergent by assumption so that we must have for
large $r$,
%
%
\begin{equation}
\PP\bigl\{\|s_n -f \| \ge6\varepsilon\mbox{ for all }n \ge
\rho^r - \rho \bigr\} = 0.
\end{equation}
It follows that
%
%
\begin{equation}
\PP\bigl\{\|s_n - f\| < 6\varepsilon\mbox{ infinitely often}\bigr\} =1,\qquad
\varepsilon>0,
\end{equation}
which implies (a).

\fbox{(a) $ \Rightarrow$ (b)} This follows directly from the
Borel--Cantelli lemma.
\end{pf}

Our next result gives a simplification of the criterion for clustering
under the additional assumption that $\{S_n/c_n\}$ is bounded in
probability, that is, we are assuming that
%
%
\begin{equation}
\forall \varepsilon>0\ \exists K_{\varepsilon} >0\qquad \PP\bigl\{|S_n|
\ge K_{\varepsilon}c_n\bigr\} < \varepsilon. \label{B}
\end{equation}
Using Theorem 1.1.5 in \cite{dlpG}, we can infer from this assumption
that also
%
%
\begin{equation}
\forall \varepsilon>0\ \exists K'_{\varepsilon} >0\qquad \PP \Bigl
\{\max_{1
\le
k \le n}|S_k| \ge K'_{\varepsilon}c_n
\Bigr\} < \varepsilon\label{B1}.
\end{equation}

\begin{propo} \label{pro35}
Under assumption $(\ref{B})$ the following are equivalent:
\begin{longlist}[(a)]
\item[(a)] $f \in C(\{s_n\dvtx n \ge1\})\ \mathrm{a.s.}$;
\item[(b)] $\sum_{n=1}^{\infty} n^{-1}\PP\{\|s_n - f\| < \varepsilon
\} =
\infty, \varepsilon>0$.
\end{longlist}
\end{propo}

\begin{pf}
\fbox{(a) $\Rightarrow$ (b)} It is obviously enough to show
that (a) implies for any $\varepsilon>0$,
%
%
\begin{equation}
\sum_{n=1}^{\infty} n^{-1}\PP\bigl\{
\|s_n -f \| < 4\varepsilon\bigl(1 + \|f\|\bigr)\bigr\} = \infty\label{ass}.
\end{equation}
Recall that by Proposition \ref{pro1} we have for \textit{any} $\rho>1$,
%
%
\begin{equation}
\sum_{k=0}^{\infty} a(\varepsilon, \rho, k) =
\infty, \label{pr}
\end{equation}
where we set
\[
a(\varepsilon, \rho, k) =\PP\bigl\{\|s_n -f \| < \varepsilon\mbox{
for some }n \in\bigl[\rho^k,\rho^{k+1}\bigr[\bigr\}.
\]
Therefore (b) follows once it has been proven that relation $(\ref
{pr})$ with a small $\rho=\rho(\varepsilon) >1$ implies $(\ref{ass})$.

To that end we first show that for $\rho^k \le m < \rho^{k+1}\le n
<\rho
^{k+2}$ and small enough $\rho>1$,
%
%
\begin{eqnarray} \label{111}
&&\bigl\{\|s_n-f\| < 4\bigl(1+\|f\|\bigr)\varepsilon\bigr\}
\nonumber
\\[-8pt]
\\[-8pt]
\nonumber
&&\qquad \supset \Bigl\{\max
_{m \le j \le n}|S_j - S_m| \le\varepsilon
c_n, \|s_m -f\| < \varepsilon\Bigr\}.
\end{eqnarray}
To verify $(\ref{111})$ observe that
\[
\|s_n -f\| = \sup_{0 \le t \le m/n}\bigl|s_n(t) -
f(t)\bigr| \vee\sup_{m/n \le
t \le1}\bigl|s_n(t) - f(t)\bigr| =:
\Delta_{n,1}^{(m)} \vee\Delta_{n,2}^{(m)}.
\]
Using the fact that $S_{(m)}(t)= S_{(n)}(mt/n), 0 \le t \le1$, it is
easy to see that
\begin{eqnarray*}
\Delta_{n,1}^{(m)}&\le&\|s_m-f\|+ \sup
_{0 \le t \le m/n}\bigl|f(t)-c_m c_n^{-1}
f(nt/m)\bigr|
\\
&\le& \|s_m-f\|+(1-c_m/c_n)\|f\| + \sup
_{0 \le t \le m/n}\bigl|f(nt/m) - f(t)\bigr|.
\end{eqnarray*}
Recall that by condition (\ref{cn2}) we have $c_m/c_n \ge(1+\varepsilon
)^{-1}m/n \ge(1+\varepsilon)^{-1}\rho^{-2}$ if $m \ge m_{\varepsilon}$.

Choose now $\rho'_{\varepsilon} > 1$ so small that $ (1+\varepsilon
)^{-1}{\rho
'_{\varepsilon}}^{-2} \ge1-2\varepsilon$.

Further, let $\delta>0$ be
small enough so that $|f(u)-f(v)| \le\varepsilon$ whenever $|u-v| <
\delta
$.

Setting $\rho_{\varepsilon} =
\rho'_{\varepsilon} \wedge(1+ \delta)^{1/2}$,
we then have if $m \ge m_{\varepsilon}$ and $1 < \rho\le\rho_{\varepsilon}$,
%
%
\begin{equation}
\|s_m -f\| < \varepsilon\Rightarrow\Delta_{n,1}^{(m)}
\le2\bigl(1+\|f\| \bigr)\varepsilon. \label{112}
\end{equation}
We now turn to the variable $ \Delta_{n,2}^{(m)}$ for which we clearly have
\begin{eqnarray*}
\Delta_{n,2}^{(m)}&\le& \sup_{m/n \le t \le
1}\bigl|c_n^{-1}
\bigl(S_{(n)}(t)-S_m\bigr)\bigr| +\bigl|S_m/c_m
- f(1)\bigr|
\\
&&{} + \sup_{m/n \le t \le1}\bigl|f(t)-(c_m/c_n)f(1)\bigr|.
\end{eqnarray*}
Arguing as above we find that
%
%
\begin{equation}\quad
\|s_m -f\| < \varepsilon\quad\Rightarrow\quad\Delta_{n,2}^{(m)}
\le\max_{m \le
j \le n}|S_j - S_m|/c_n
+ 2\varepsilon\bigl(1 + \|f\|\bigr), \label{113}
\end{equation}
provided that $m \ge m_{\varepsilon}$ and
$\rho\le\rho_{\varepsilon}$.\vadjust{\goodbreak}

Combining $(\ref{112})$ and $(\ref{113})$ we get $(\ref{111})$.

Let $\tau_k$ and $I_k$ be defined as in the proof of Proposition \ref
{pro1}. Then we have for large enough $k$,
\begin{eqnarray*}
&&\sum_{n \in I_{k+1}} \PP\bigl\{\|s_n -f \| <
4\bigl(1+\|f\|\bigr)\varepsilon\bigr\}
\\
&&\qquad\ge\sum_{m \in I_k} \sum_{n \in I_{k+1}}
\PP\bigl\{\|s_n -f\| < 4\bigl(1+\| f\| \bigr)\varepsilon, \tau_k =m
\bigr\}
\\
&&\qquad\ge\sum_{m \in I_k} \sum_{n \in I_{k+1}}
\PP\Bigl\{\max_{m \le j \le
n}|S_j -S_m| \le
\varepsilon c_n\Bigr\}\PP\{\tau_k =m\}
\\
&& \qquad\ge\bigl\{\rho^k (\rho-1) -1\bigr\}\PP\Bigl\{\max
_{1 \le j \le r_k} |S_j| \le \varepsilon c_{n_k}
\Bigr\}a(\varepsilon,\rho,k),
\end{eqnarray*}
where $r_k \le\rho^{k+2}- \rho^k +2$ and $n_k \ge\rho^{k+1}-1$.
Noticing that $c_{n_k}/c_{r_k} \ge(n_k/r_k)^{1/2}$, with $ \liminf_{k
\rightarrow\infty} (n_k/r_k)^{1/2} \ge\rho^{1/2}/(\rho^2-1)^{1/2}$,
and recalling (\ref{B1}), we can choose a constant $1 < \bar{\rho
}_{\varepsilon} < \rho_{\varepsilon}$ so that we have for $1 < \rho< \bar
{\rho}_{\varepsilon}$ and large $k$,
\[
\PP\Bigl\{\max_{1 \le j \le r_k}|S_j| \le\varepsilon
c_{n_k}\Bigr\} \ge \PP\Bigl\{\max_{1 \le j \le r_k}
|S_j| \le K c_{r_k}\Bigr\}\ge1/2.
\]
Consequently, we have for large $k$ and $1 < \rho< \bar{\rho
}_{\varepsilon}$,
\[
\sum_{n \in I_{k+1}}n^{-1} \PP\bigl\{
\|s_n -f\| < 4\bigl(1+\|f\|\bigr)\varepsilon\bigr\} \ge \frac{1}{2\rho^2}\bigl(
\rho-1-\rho^{-k}\bigr)a(\varepsilon,\rho,k),
\]
which implies $(\ref{ass})$ and thus (b).

\fbox{(b) $\Rightarrow$ (a)} Noting that we have for any $\rho>1$,
\[
\sum_{n \in I_k} n^{-1} \PP\bigl\{\|s_n
-f\| < \varepsilon\bigr\} \le\bigl(\rho- 1 + \rho^{-k}\bigr)\PP\bigl\{
\|s_n - f\| < \varepsilon\mbox{ for some }n \in I_k\bigr\},
\]
this implication follows immediately from Proposition \ref{pro1}.
\end{pf}
\section{Clustering in $\R^d$}\label{sec4}

In this section we look at $d$-dimensional random vectors, where again
$|\cdot|$ will denote the Euclidean norm. We first provide a criterion
for clustering in the functional case in terms of Brownian motion. We
use the following strong approximation result.

\begin{theorem}[(Sakhanenko \cite{sakh})]\label{thmB}
Let $X_1^*, \ldots, X_n^*$ be
independent mean zero random vectors in $\R^d$ and assume that $\E
|X^*_i|^p < \infty, 1 \le i \le n$ for some $p \in\,]2,3]$. Let $x > 0$
be fixed. If the underlying probability space is rich enough, one can
construct independent normally distributed mean zero\vadjust{\goodbreak} random vectors
$Y^*_1, \ldots, Y^*_n$ with $\operatorname{cov} (X^*_i) =\operatorname{cov}(Y^*_i),
1 \le i \le n$ such that
\[
\PP \Biggl\{\max_{1 \le k \le n} \Biggl\llvert \sum
_{j=1}^k \bigl(X_j^* -
Y^*_j\bigr)\Biggr\rrvert \ge x \Biggr\} \le K \sum
_{i=1}^n \E\bigl|X^*_i\bigr|^p
/x^p,
\]
where $K$ is a positive constant depending on $d$ only.
\end{theorem}

Note that there is no assumption on the covariance matrices of the
random vectors $X_1^*, \ldots, X_n^*$. This will be crucial for the
subsequent proof since we will apply it to truncated random vectors
where the original (``untruncated'') random vectors do not need to have
finite covariance matrices.

In this way we obtain the following criterion for clustering in the
functional~LIL:
%
\begin{theo} \label{cri} Let $X=(X^{(1)},\ldots,X^{(d)}) \dvtx \Omega\to
\R
^d$ be a mean zero random vector, and let $\{c_n\}$ be a sequence of
positive real numbers satisfying conditions~(\ref{cn1}) and (\ref{cn2}).
Set $s_n = S_{(n)}/c_n\dvtx \Omega\to C_d[0,1]$. Assuming that $\sum_{n=1}^{ \infty}\PP\{|X| \ge c_n\} <\infty$, the following are equivalent:
\begin{longlist}[(a)]
\item[(a)] $f \in C(\{s_n\dvtx n \ge1\})\ \mathrm{a.s.}$;
\item[(b)] we have for any $\varepsilon> 0$,
\[
\sum_{n=1}^{\infty} n^{-1}\PP\bigl\{\|
\Gamma_n W_{(n)}/c_n - f\| < \varepsilon \bigr\} =
\infty,
\]
where $\Gamma_n$ is the positive semidefinite symmetric matrix such
that
\[
\Gamma_n^2 = \bigl(\E\bigl[X^{(i)}X^{(j)}
I\bigl\{|X| \le c_n\bigr\}\bigr]\bigr)_{1 \le i, j \le d}
\]
and $W_{(n)}(t)= W(nt), 0 \le t \le1$ with $W$ being a standard
$d$-dimensional Brownian motion.
\end{longlist}
\end{theo}
In the proof we make extensive use of the following lemma. The easy
proof of this lemma is omitted.
%
\begin{lem} \label{le41}
Let $\xi_n, \eta_n\dvtx \Omega\to C_d [0,1]$ be random elements such that
\[
\sum_{n=1}^{\infty} n^{-1}\PP\bigl\{\|
\xi_n -\eta_n\| \ge\varepsilon\bigr\} < \infty, \qquad\varepsilon>0.
\]
Then we have for any function $f \in C_d [0,1]$,
\[
\sum_{n=1}^{\infty} n^{-1}\PP\bigl\{\|
\xi_n -f\| < \varepsilon\bigr\} < \infty\qquad \forall\varepsilon>0
\]
if and only if
\[
 \sum_{n=1}^{\infty}
n^{-1}\PP\{\|\eta_n -f\| < \varepsilon\} < \infty\qquad\forall
\varepsilon>0.
\]
\end{lem}
We record the following facts which can be proved similarly as in the
1-dimensional case (refer to Lemma 1 in \cite{EL-1}).

If $X$ is a mean zero random vector such that $\sum_{n=1}^{\infty}
\PP\{
|X| \ge c_n\}< \infty$, where $c_n$ satisfies the two conditions (\ref
{cn1}) and (\ref{cn2}), we have:

\textit{Fact} 1. $\sum_{n=1}^{\infty} \E[|X|^3 I\{|X| \le c_n\}]/c_n^3 <
\infty$;

\textit{Fact} 2. $\E[|X|I\{|X|\ge c_n\}] = o(c_n/n)$ as $n \to\infty$;

\textit{Fact} 3. $\E[|X|^2 I\{|X| \le c_n\}] = o(c_n^2 /n)$ as $n \to
\infty
$.

We are ready to prove Theorem \ref{cri}. By a slight abuse of notation
we also denote the Euclidean matrix norm by $\|\cdot\|$ if $\Gamma$ is
a $(d,d)$-matrix. That is, we set $\|\Gamma\|=\sup_{|x|\le1}|\Gamma
x|$. Recall that if $\Gamma$ is a symmetric matrix, $\|\Gamma\|^2$ is
equal to the largest eigenvalue of the matrix $\Gamma^2$.

\begin{pf*}{Proof of Theorem \ref{cri}}
(i) Set $X'_{n, j} = X_j I\{|X_j| \le c_n\}, X^*_{n,j} =
X'_{n,j} - \E X'_{n,j}, 1 \le j \le n, n \ge1$, and let $S^*_{(n)}$ be
the partial sum process based on $X^*_{n,1},\ldots,X^*_{n,n} $. Finally
set $s^*_n = S^*_{(n)}/c_n, n \ge1$.

Then we have
%
%
\begin{equation}
\label{21} \sum_{n=1}^{\infty} n^{-1}
\PP\bigl\{\bigl\|s_n - s^*_n\bigr\| \ge\varepsilon\bigr\} < \infty,\qquad
\varepsilon> 0.
\end{equation}
To verify (\ref{21}) observe that
\begin{eqnarray*}
\bigl\|s_n - s^*_n\bigr\| &\le& \max_{1 \le k \le n}
\Biggl\llvert \sum_{j=1}^k
\bigl(X_j - X^*_j\bigr)\Biggr\rrvert\Big /c_n
\\
&\le& \Biggl(\sum_{j=1}^n
|X_j| I\bigl\{|X_j| > c_n\bigr\} + n\E|X|I
\bigl\{|X|>c_n\bigr\} \Biggr)\big/c_n.
\end{eqnarray*}
Recalling Fact 2 we get for large $n$,
%
%
\begin{eqnarray}
\label{23} \PP\bigl\{\bigl\|s_n - s^*_n\bigr\| \ge\varepsilon\bigr
\} &\le& \PP \Biggl(\sum_{j=1}^n
|X_j| I\bigl\{|X_j| > c_n\bigr\} >
\frac{\varepsilon}{2}c_n \Biggr)
\nonumber
\\[-8pt]
\\[-8pt]
\nonumber
& \le& n\PP\bigl\{|X| \ge c_n\bigr\},
\end{eqnarray}
and we see that (\ref{21}) holds.

Noting that Facts 2 and 3 also imply that $\E|S_n| /c_n \to0$ (see
Lemma 1, \cite{EL}), we trivially have that $\{S_n/c_n\dvtx n \ge1\}$ is
stochastically bounded. Consequently, Proposition \ref{pro35} can be
applied which in combination with (\ref{21}) and Lemma \ref{le41} gives
%
%
\begin{eqnarray}
\label{22} f \in C\bigl(\{s_n\dvtx n \ge1\}\bigr) \quad\mbox{a.s.}\quad
\Longleftrightarrow\quad
\sum_{n=1}^{\infty}
n^{-1}\PP\bigl\{\bigl\|s^*_n - f\bigr\| < \varepsilon\bigr\} =
\infty, \quad\varepsilon>0.
\nonumber
\\[-8pt]
\end{eqnarray}

(ii) In this part we will use Theorem \ref{thmB}. From Fact 1 it easily follows
that one can find a sequence $\tilde{c}_n \nearrow\infty$ so that
$\tilde{c}_n /c_n \to0$ as $n \to\infty$ and we still have
%
%
\begin{equation}
\label{2.4} \sum_{n=1}^{\infty} \E
\bigl[|X|^3 I\bigl\{|X| \le c_n\bigr\}\bigr]/\tilde{c}_n^3
< \infty.
\end{equation}
Let $ n \ge1$ be fixed. Employing the afore-mentioned result along
with the $c_r$-inequality, we can construct independent
$N(0,I)$-distributed random vectors $Y_{n,1}, \ldots, Y_{n,n}$ such
that we have
%
%
\begin{equation}
\label{2.5} \PP \Biggl\{\max_{1 \le k \le n} \Biggl\llvert \sum
_{j=1}^k \bigl(X^*_j -
\Gamma^*_n Y_{n,j}\bigr)\Biggr\rrvert \ge
\tilde{c}_n \Biggr\} \le8K n \E\bigl[|X|^3I\bigl\{|X| \le
c_n\bigr\}\bigr]/\tilde{c}_n^3,
\end{equation}
where $\Gamma^*_n$ is the symmetric positive semidefinite matrix such
that $(\Gamma^*_n)^2 = \operatorname{cov}(X^*_{n,1})$.

Letting $T_{(n)}$ be the partial sum process based on the random
vectors $Y_{n,1},\ldots, Y_{n,n}$, $t_n = T_{(n)}/c_n$ and recalling
(\ref{2.4}), we find that
%
%
\begin{equation}
\sum_{n=1}^{\infty} n^{-1}
\PP\bigl\{\bigl\|s^*_n - \Gamma^*_n t_n \bigr\| \ge
\varepsilon c_n\bigr\} <\infty,\qquad \varepsilon>0.
\end{equation}
This means in view of Lemma \ref{le41} and relation (\ref{22}) that
%
%
\begin{eqnarray}
\label{26} f \in C\bigl(\{s_n\dvtx n \ge1\}\bigr)\quad \mbox{a.s.}\quad
\Longleftrightarrow\quad
\sum_{n=1}^{\infty}
n^{-1}\PP\bigl\{\bigl\|\Gamma^*_n t_n - f\bigr\| <
\varepsilon\bigr\} = \infty, \quad\varepsilon>0.
\nonumber
\\[-8pt]
\end{eqnarray}

(iii) Let $W^{[n]}(t), t \ge0$ be a Brownian motion satisfying
\[
W^{[n]}(k) = \sum_{j=1}^k
Y_{n,j},\qquad 1 \le k \le n.
\]
Then we have
\[
\bigl\|T_{(n)} - W^{[n]}_{(n)}\bigr\| \le2 \max
_{0 \le j \le n-1} \sup_{0 \le u
\le1} \bigl|W^{[n]}(j+u)-W^{[n]}(j)\bigr|,
\]
and we can conclude for $x >0$,
\[
\PP\bigl\{\bigl\|T_{(n)} - W^{[n]}_{(n)}\bigr\| \ge x\bigr\}
\le n\PP \Bigl\{\sup_{0 \le u
\le1}\bigl|W^{[n]}(u)\bigr| \ge x/2 \Bigr\}
\le2dn\exp \biggl(-\frac
{x^2}{8d} \biggr).
\]
It follows that
\[
\PP\bigl\{\bigl\|\Gamma_n^* \bigl(t_n - W^{[n]}_{(n)}/c_n
\bigr)\bigr\| \ge\varepsilon\bigr\} \le 2dn\exp \biggl(-\frac{\varepsilon^2 c_n^2}{8d\|\Gamma^*_n\|^2} \biggr).
\]
As $\|\Gamma^*_n\|^2 \le\E|X^*_{n,1}|^2 \le\E[|X|^2 I\{|X| \le c_n\}]
= o(c_n^2/n)$ (see Fact 3), we readily obtain that
\[
\sum_{n=1}^{\infty} n^{-1} \PP\bigl\{
\bigl\|\Gamma_n^* \bigl(t_n - W^{[n]}_{(n)}/c_n
\bigr)\bigr\| \ge\varepsilon\bigr\} < \infty,\qquad \varepsilon>0.
\]
Consequently we have by Lemma \ref{le41} and (\ref{26}) and since
$W^{[n]}_{(n)} \stackrel{d}{=}W_{(n)}, n \ge1$,
%
%
\begin{eqnarray}
\label{27} f \in C\bigl(\{s_n\dvtx n \ge1\}\bigr)\quad \mbox{a.s.}\quad\Longleftrightarrow\quad
\sum_{n=1}^{\infty}
n^{-1}\PP\bigl\{\bigl\|\Gamma^*_n w_n - f\bigr\| <
\varepsilon\bigr\} = \infty, \quad\varepsilon>0,
\nonumber
\\[-8pt]
\end{eqnarray}
where $w_n = W_{(n)}/c_n$.

(iv) Observing that $\Delta_n = \Gamma_n^2 - \Gamma_n^{*2}$ is a
positive semidefinite symmetric matrix, we clearly have
\[
\Gamma_n W_{(n)} \stackrel{d} {=} \Gamma^*_n
W_{(n)} + \bar{\Delta}_n \bar{W}_{(n)}=:Z_n,
\]
provided that $\bar{W}_{(n)}(t)=\bar{W}(nt), 0 \le t \le1$, where
$\bar
{W}(s), s \ge0$ is another Brownian motion which is independent of
$W$, and $\bar{\Delta}_n$ is the positive semidefinite symmetric matrix
satisfying $ \bar{\Delta}_n^2 = \Delta_n$.

It follows that
\[
\PP\bigl\{\bigl\|Z_n - \Gamma^*_n W_{(n)} \bigr\| \ge
\varepsilon c_n\bigr\} \le\PP\bigl\{\| \bar {\Delta}_n\| \|
\bar{W}_{(n)}\| \ge\varepsilon c_n\bigr\}.
\]
Since we have $\bar{W}_{(n)}(t) \stackrel{d}{=}\sqrt{n}W(t), 0 \le t
\le1$, we find that this probability is
\[
\le\PP \Bigl\{\sup_{0 \le t \le1}\bigl|W(t)\bigr| \ge\varepsilon c_n/
\bigl(\sqrt {n}\| \bar{\Delta}_n\| \bigr) \Bigr\} \le2d\exp \biggl(-
\frac{\varepsilon^2 c_n^2}{2dn\|\bar{\Delta}_n\|
^2} \biggr).
\]
By the definition of the matrix norm we further have
\[
\|\bar{\Delta}_n\|^2 = \mbox{ largest eigenvalue of }
\Delta_n = \sup_{|t|=1} \langle t,
\Delta_n t\rangle.
\]
A straightforward calculation gives if $|t| \le1$,
\begin{eqnarray*}
\langle t, \Delta_n t\rangle&=& \Biggl(\E \sum
_{i=1}^d t_i X^{(i)} I\bigl\{ |X|
\le c_n\bigr\} \Biggr)^2 = \Biggl(\E \sum
_{i=1}^d t_i X^{(i)} I\bigl\{|X| >
c_n\bigr\} \Biggr)^2
\\
&\le& \Biggl( \sum_{i=1}^d
|t_i| \Biggr)^2 \bigl(\E|X| I\bigl\{|X| > c_n\bigr\}
\bigr)^2 \le d \bigl(\E|X| I\bigl\{|X| > c_n\bigr\}
\bigr)^2.
\end{eqnarray*}
Recalling Fact 2 we see that $\|\bar{\Delta}_n\|^2 = o(c_n^2/ n^2)$ as
$n \to\infty$, which in turn implies that
%
%
\begin{equation}
\sum_{n=1}^{\infty}n^{-1} \PP\bigl\{
\bigl\|\Gamma^*_n w_n - Z_n/c_n\bigr\| \ge
\varepsilon \bigr\} < \infty,\qquad \varepsilon>0.
\end{equation}
Using once more Lemma \ref{le41} along with the fact that $Z_n/c_n
\stackrel{d}{=}\Gamma_n w_n$, we get that
%
%
\begin{eqnarray}
f \in C\bigl(\{s_n\dvtx n \ge1\}\bigr)\quad \mbox{a.s.}\quad \Longleftrightarrow\quad
\sum_{n=1}^{\infty} n^{-1}\PP\bigl\{\|
\Gamma_n w_n - f\| < \varepsilon\bigr\} = \infty,\quad
\varepsilon>0,
\nonumber
\\[-8pt]
\end{eqnarray}
and Theorem \ref{cri} has been proven.
\end{pf*}

We next look at the case where the random vector $X\dvtx \Omega\to\R^d$
has independent components. In this case we can prove the following:
%
\begin{theo} \label{crimod}
Let $X=(X^{(1)},\ldots,X^{(d)})\dvtx \Omega\to\R^d$ be a mean zero
random vector such that $X^{(1)},\ldots,X^{(d)}$ are independent.

Assuming that $\sum_{n=1}^{ \infty}\PP\{|X| \ge c_n\} <\infty$, where
$c_n$ is as in (\ref{cri}), the following are equivalent:
\begin{longlist}[(a)]
\item[(a)] $f=(f_1,\ldots,f_d) \in C(\{s_n\dvtx n \ge1\})\ \mathrm{a.s.}$;
\item[(b)] we have for any $\varepsilon> 0$,
\[
\sum_{n=1}^{\infty} n^{-1}\prod
_{i=1}^d \PP\bigl\{\bigl\|f_i -
\sigma_{n,i} W'_{(n)}/c_n \bigr\| <
\varepsilon\bigr\} = \infty,
\]
where $\sigma_{n,i}^2 = \E[(X^{(i)})^2 I\{|X^{(i)}| \le c_n\}], 1 \le i
\le d$ and $W'_{(n)}(t)= W'(nt), 0 \le t \le1$ with $W'$ being a
standard 1-dimensional Brownian motion.
\end{longlist}
\end{theo}
The proof is similar to the previous one and we will just indicate the
changes.

\begin{pf*}{Proof of Theorem \ref{crimod}}
(i) We define the random vectors $X'_{n,j}, 1 \le j \le n$
as follows:
\[
X'_{n,j} = \bigl(X^{(1)} I\bigl
\{\bigl|X^{(1)}\bigr| \le c_n\bigr\},\ldots, X^{(d)} I\bigl
\{\bigl|X^{(d)}\bigr| \le c_n\bigr\}\bigr), \qquad 1 \le j \le n, n \ge1.
\]
Letting again $X^*_{n,j} = X'_{n,j} - \E X'_{n,j}, 1 \le j \le n, n \ge
1$, we have
\begin{eqnarray*}
\bigl\|s_n - s^*_n\bigr\|& =& \max_{1 \le k \le n} \Biggl
\llvert \sum_{j=1}^k \bigl(X_j
- X^*_j\bigr)\Biggr\rrvert\Big/c_n
\\
&\le& \sum_{i=1}^d \Biggl(\sum
_{j=1}^n \bigl|X^{(i)}_j\bigr| I\bigl
\{\bigl|X^{(i)}_j\bigr| > c_n\bigr\} + n
\E\bigl|X^{(i)}\bigr|I\bigl\{\bigl|X^{(i)}\bigr|>c_n\bigr\}
\Biggr)\Big/c_n
\\
&\le& d \Biggl(\sum_{j=1}^n
|X_j| I\bigl\{|X_j| > c_n\bigr\} + n\E|X|I\bigl\{
|X|>c_n\bigr\} \Biggr)\Big/c_n,
\end{eqnarray*}
and as in the previous proof we see that we can replace $s_n$ by
$s^*_n$.

(ii) This part remains essentially unchanged. Note that $\Gamma^*_n$ is
now a diagonal matrix. The only difference is that we have to use a
slightly different upper bound for $\E|X^*_{n,1}|^3$,
\[
\E\bigl|X^*_{n,1}\bigr|^3 \le8\E\bigl|X'_{n,1}\bigr|^3
\le8d^{1/6}\sum_{i=1}^d \E
\bigl|X^{(i)}\bigr|^3 I\bigl\{\bigl|X^{(i)}\bigr| \le c_n
\bigr\},
\]
where the second bound easily follows from the H\"older inequality.

Applying Fact 1 for each $X^{(i)}$ we see that $\sum_{n=1}^{\infty}
\E|
X'_{n,1}|^3/c_n^3 <\infty$.

(iii) Here we use the fact\vspace*{1pt} that $\|\Gamma_n^*\|^2 \le\max_{1 \le i
\le
d} \E(X^{(i)})^2 I\{|X^{(i)}| \le c_n\} = o(c_n^2/n)$ since we can
employ Fact 3\vspace*{1pt} for the (finitely many) random variables $X^{(i)}, 1 \le
i \le d$ as well.

(iv) Since $\Delta_n$ is a diagonal matrix, we have that
\[
\|\Delta_n\| =\max_{1 \le i \le d} \bigl(
\E\bigl[X^{(i)} I\bigl\{\bigl|X^{(i)}\bigr| \le c_n\bigr\}
\bigr]\bigr)^2 = \max_{1 \le i \le d} \bigl(\E\bigl[X^{(i)} I
\bigl\{\bigl|X^{(i)}\bigr| > c_n\bigr\}\bigr]\bigr)^2,
\]
which is of order $o(c_n^2/n^2)$ due to Fact 2 [applied for the
components $X^{(i)}$].

(v) (small extra step). We have shown so far that
%
%
\begin{eqnarray}
f &\in& C\bigl(\{s_n\dvtx n \ge1\}\bigr)\quad \mbox{a.s.}\quad
\Longleftrightarrow\quad
\sum_{n=1}^{\infty} n^{-1}\PP\bigl\{
\|D_n w_n - f\| < \varepsilon\bigr\} = \infty,\quad \varepsilon>0,
\nonumber
\\[-8pt]
\end{eqnarray}
where $D_n=\Gamma_n=\operatorname{diag}(\sigma_{n,1},\ldots,\sigma
_{n,d})$ is
a diagonal matrix.

It is trivial that we can replace the sup-norm $\|\cdot\|$ based on the
Euclidean norm $|\cdot|$ by the equivalent sup-norm $\|\cdot\|_{+}$
which is based on the norm $|x|_{+} = \max_{1 \le i \le d} |x_i|$. In
this case we also have for $g=(g_1,\ldots,g_d)$ that $\|g\|_{+} = \max_{1 \le i \le d}\sup_{0 \le t \le1}|g_i(t)|$.
Thus we have
\begin{eqnarray*}
&&f \in C\bigl(\{s_n\}\bigr)\qquad \mbox{a.s.} \\
&&\qquad\Longleftrightarrow \quad\sum
_{n=1}^{\infty} n^{-1}\PP\bigl\{
\Vert D_n w_n - f\Vert _{+} < \varepsilon\bigr\} = \infty,\qquad
\varepsilon>0
\\
&&\qquad\Longleftrightarrow\quad \sum_{n=1}^{\infty}
n^{-1}\PP \Biggl(\bigcap_{i=1}^d
\bigl\{\bigl\| \sigma_{n,i} w^{(i)}_n - f_i
\bigr\| < \varepsilon\bigr\} \Biggr) = \infty, \qquad\varepsilon>0,
\end{eqnarray*}
and Theorem \ref{crimod} follows by independence.\vadjust{\goodbreak}
\end{pf*}

Analogous results hold for the cluster sets $A=C(\{S_n/c_n\dvtx n \ge1\})$.
%
\begin{theo} \label{criA} Let $X\dvtx \Omega\to\R^d$ be a mean zero
random vector, and let $\{c_n\}$ be a sequence of positive real numbers
satisfying conditions (\ref{cn1}) and (\ref{cn2}).
Assuming that $\sum_{n=1}^{ \infty}\PP\{|X| \ge c_n\} <\infty$, the
following are equivalent:
\begin{longlist}[(a)]
\item[(a)] $x \in C(\{S_n/c_n\dvtx n \ge1\})$ a.s.;
\item[(b)] we have for any $\varepsilon> 0$,
\[
\sum_{n=1}^{\infty} n^{-1}\PP\bigl\{\bigl|
\Gamma_n W(n)/c_n - x\bigr| < \varepsilon \bigr\} = \infty,
\]
where $\Gamma_n$ is as in Theorem \ref{cri}, and $W$ is a standard
$d$-dimensional Brownian motion.
\end{longlist}
Furthermore, if $X$ has independent components $X^{(1)}, \ldots,
X^{(d)}$, \textup{(a)} is also equivalent to the following:
\begin{longlist}[(c)]
\item[(c)] we have for any $\varepsilon> 0$
\[
\sum_{n=1}^{\infty} n^{-1}\prod
_{i=1}^d \PP\bigl\{\bigl|x_i -
\sigma_{n,i} W'(n)/c_n \bigr| < \varepsilon\bigr\} =
\infty,
\]
where $\sigma_{n,i}^2, 1 \le i \le d$ is as in Theorem \ref{crimod},
and $W'$ is a standard 1-dimensional Brownian motion.
\end{longlist}
\end{theo}

\begin{pf} Using a version of Lemma \ref{le41} for random vectors and
recalling relation (\ref{Kes}), the equivalence of (a) and (b) follows
once it has been shown that
%
%
\begin{equation}
\label{411} \sum_{n=1}^{\infty} n^{-1}
\PP\bigl\{\bigl|S_n - \Gamma_n W(n)\bigr| \ge\varepsilon
c_n\bigr\} < \infty, \qquad\varepsilon>0.
\end{equation}
From the proof of Theorem \ref{cri} it follows that we actually have
\[
\sum_{n=1}^{\infty} n^{-1} \PP\bigl\{
\|S_{(n)} - \Gamma_n W_{(n)}\| \ge \varepsilon
c_n\bigr\} < \infty,\qquad \varepsilon>0,
\]
which trivially implies (\ref{411}).

The proof of the equivalence of (a) and (c) is similar.
\end{pf}
\section{\texorpdfstring{Proof of Theorem \protect\ref{theo21}}
{Proof of Theorem 2.1}}\label{sec5}
Let $S_n^{(i)}$ and
$S_{(n)}^{(i)}$ denote the $i$th coordinate of $S_n$ and $S_{(n)}$,
respectively. Note that $S_{(n)}^{(i)}$ is then the 1-dimensional
partial sum process based on the sequence $S_n^{(i)}, n \ge1$.

From Theorem \ref{thmA} it follows that $\limsup_{n \to\infty}|S_n|/c_n =
\alpha
_0 < \infty$ with probability one, which clearly implies that $A=C(\{
S_n/c_n\})$ is a compact subset of $\R^d$. Applying Proposition \ref
{pro33} we then have $\mathcal{A}$ compact in $C_d[0,1]$ and also that
both $A$ and $\mathcal{A}$ are nonempty.

Furthermore, $\alpha_{i} \le\alpha_0 <\infty$ for $i=1,\ldots,d$,
whence by Theorem 3 of \cite{EL-1} with probability one,
$\limsup_{n \rightarrow\infty}S_n^{(i)}/c_n = \alpha_i, i=1,\ldots,d$.
This in turn implies by Theorem 1 of \cite{E} that with probability one
\[
C\bigl(\bigl\{S_{(n)}^{(i)}/c_n\dvtx n \ge1\bigr\}
\bigr)=\alpha_{i}\mathcal{K}
\]
and
\[
\lim_{n \rightarrow\infty} \inf_{h_i \in\alpha_{i} \mathcal{K}} \bigl\Vert
h_i - S_{(n)}^{(i)}/c_n\bigr\Vert =0,\qquad 1 \le
i \le d.
\]
Therefore,
with probability one
\[
\lim_{n \rightarrow\infty} \sum_{i=1}^d
\inf_{h_i \in\alpha_{i}
\mathcal{K}} \bigl\Vert h_i - S_{(n)}^{(i)}/c_n
\bigr\Vert =0.
\]
Using $\Vert f-g\Vert  \le\sum_{i=1}^d \Vert f_i-g_i\Vert $ for $f=(f_1,\ldots
,f_d),g=(g_1,\ldots,g_d) \in C_d([0,1])$, we have
\[
\limsup_{n \rightarrow\infty} \inf_{h \in\alpha_{1}\mathcal{K}
\times\cdots\times\alpha_{d}\mathcal{K}}
\|h-S_{(n)}/c_n\| \le \lim_{n \rightarrow\infty}\sum
_{i=1}^d \inf_{h_i \in\alpha_{i}
\mathcal
{K}}
\bigl\Vert h_i - S_{(n)}^{(i)}/c_n\bigr\Vert =0,
\]
and we see that (\ref{24}) holds since $ \alpha_{1}\mathcal{K}
\times
\cdots\times\alpha_{d}\mathcal{K}$ is a compact subset of
$C_d([0,1])$.

To prove the other inclusion in Theorem \ref{theo21} we need more
notation. As the matrices $\Gamma_n$ defined in Theorem \ref{cri} are
positive semidefinite and symmetric, we can find orthonormal bases $\{
u_{n,1}, \ldots, u_{n,d}\}$ of $\R^d$ consisting of eigenvectors of
$\Gamma_n$. Let $\lambda_{n,i}$ be the corresponding eigenvalues. We
can assume w.l.o.g. that $\lambda_{n,i}, 1 \le i \le d_n$ are the
nonzero eigenvalues, where $1 \le d_n \le d$.
Set $\xi_n = XI\{|X| \le c_n\}, n \ge1$. Note that then by definition
of $\Gamma_n^2$, $\E[\langle\xi_n, u_{n,i}\rangle^2]=\langle u_{n,i},
\Gamma_n^2 u_{n,i}\rangle= \lambda_{n,i}^2=0, d_n < i \le d$.

Thus with probability one,
$\xi_n=\sum_{i=1}^d \langle\xi_n, u_{n,i}\rangle u_{n,i} = \sum_{i=1}^{d_n} \langle\xi_n, u_{n,i}\rangle u_{n,i}$.

We see that $\PP\{\xi_n \in V_n\} =1$ if $V_n$ is the $d_n$-dimensional
subspace of $\R^d$ spanned by
$u_{n,i}, 1 \le i \le d_n$.

Further note that the sequence $\Gamma_n^2$ is monotone; that is,
$\Gamma_n^2 - \Gamma_m^2$ is positive semidefinite if $n \ge m$. Let
$V'_n$ be the vector space spanned by $u_{n,i}, d_n < i \le d$ if $d_n
< d$ and $\{0\}$ otherwise. This is the zero space of the quadratic
form determined by $\Gamma_n^2$, and thus by
monotonicity of $\Gamma_n^2$ we get that $V'_1 \supset V'_2 \supset
\cdots.$ As $V'_n$ is the orthogonal complement of $V_n$ we can conclude
that $V_1 \subset V_2 \subset\cdots.$ Thus there are at most $d+1$
different vector spaces [with $0 \le\operatorname{dim}(V_{n_i}) \le d$] in
this sequence, and we have $V_n = V$ eventually for some subspace $V$
of $\R^d$ with dimension $1 \le d' \le d$. Notice also that $X$ is
supported by this vector space as we have $\PP\{X \in V\}
=\lim_{n \to\infty} \PP\{\xi_n \in V\}=1$.

We first infer from Theorem \ref{cri} the following lemma.
%
\begin{lem} \label{lem50}
Under the assumptions of Theorem \ref{theo21} we have
$f=(f_1,\ldots,\break  f_d) \in\mathcal{A}$ if and only if
%
%
\begin{equation}
\label{51} \sum_{n \in\mathbb{N}_{\varepsilon}} n^{-1} \exp \Biggl(-
\sum_{i=1}^{d_n}\frac{I(\langle u_{n,i}, f\rangle_{\varepsilon})
c_n^2}{2n\lambda_{n,i}^2} \Biggr) =
\infty \qquad\forall\varepsilon>0,
\end{equation}
where $\mathbb{N}_{\varepsilon}=\{n \ge1\dvtx \|\langle u_{n,i}, f\rangle
\| <
\varepsilon, i >d_n\}$.
\end{lem}
\begin{pf}
If $U_n$ denotes the orthogonal matrix whose $i$th column is the $i$th
eigenvector $u_{n,i}$, then since the probability law of $W$ is the
same as that of $U_nW$, we have
\[
\PP\bigl(\|\Gamma_nW_{(n)}/c_n - f\| < \varepsilon\bigr)=
\PP\bigl(\|\Gamma _nU_nW_{(n)}/c_n - f\|
< \varepsilon\bigr).
\]
In addition, since the transposed matrix $U_n'$ is orthogonal, it
preserves distances given by the Euclidean norm and hence
\[
P\bigl(\|\Gamma_nU_nW_{(n)}/c_n - f\| <
\varepsilon\bigr)=\PP\bigl(\bigl\|U_n'\Gamma _nU_nW_{(n)}/c_n
-U_n' f\bigr\| < \varepsilon\bigr).
\]
Note that $D_n=U_n'\Gamma_nU_n$ is a diagonal matrix whose $i$th
diagonal entry is the eigenvalue $\lambda_{n,i}$.
Replacing the
sup-norm $\|\cdot\|$ in Theorem \ref{cri} by the equivalent norm
\[
\|g\|_{+} = \max_{1 \le i \le d} \|g_i\|,\qquad g
\in C_d[0,1],
\]
we can infer that $f \in\mathcal{A}$ if and only if
%
%
\begin{equation}
\sum_{n=1}^{\infty} n^{-1} \PP\bigl(
\bigl\|D_nW_{(n)}/c_n -U_n'
f\bigr\|_{+} < \varepsilon \bigr)=\infty, \qquad\varepsilon>0,
\end{equation}
which by independence is the same as
%
%
\begin{equation}
\label{53} \sum_{n \in\mathbb{N}_{\varepsilon}}n^{-1} \prod
_{i=1}^{d_n} \PP\bigl(\bigl\| \lambda _{n,i}
W^{(i)}_{(n)}/c_n - \langle u_{n,i}, f
\rangle\bigr\| < \varepsilon \bigr)=\infty,\qquad \varepsilon>0,
\end{equation}
where $W_{(n)}^{(i)}$ is the $i$th coordinate of $W_{(n)}$.

Now eventually in $n$ we have for $1 \le i \le d_n$,
\begin{eqnarray*}
\PP\bigl(\bigl\Vert \langle u_{n,i}, f\rangle- \lambda_{n,i}W^{(i)}_{(n)}/c_n
\bigr\Vert \leq2\varepsilon\bigr)&\geq& \PP\bigl(\bigl\Vert \langle u_{n,i}, f
\rangle_{\varepsilon} - \lambda_{n,i}W^{(i)}_{(n)}/c_n
\bigr\Vert \leq\varepsilon\bigr)
\\
&\geq& \frac{1}{2}\exp \biggl(-I\bigl( \langle u_{n,i}, f\rangle
_{\varepsilon
}\bigr)\frac{ c_n^2}{2n\lambda_{n,i}^2} \biggr).
\end{eqnarray*}
The second inequality above follows from (4.16) in Theorem 2 of \cite
{KLL} with $\alpha=0$ where we use the fact that $\lim_{n \rightarrow
\infty} c_n/(\sqrt n \lambda_{n,i}) =\infty$ for $1 \le i \le d_n$.

This last statement is true since $\lim_{n \rightarrow\infty}
c_n/(\sqrt n \sigma_{n,i}) =\infty$ for $1 \le i \le d$, which follows
from (3.3) in Lemma 1 of \cite{EL-1}. Since $\sum_{i=1}^d \lambda
_{n,i}^2=\sum_{i=1}^d \sigma_{n,i}^2$ we have $\max_{1 \leq i \leq d}
\lambda_{n,i} \leq d^{{1}/{2}}\max_{1 \leq i \leq d} \sigma_{n,i}$
whence $ c_n/(\sqrt n \lambda_{n,i}) \to\infty$ for $1 \le i \le
d_n$.

We also have from (4.17) of \cite{KLL}, with $\alpha=1$ and
$i=1,\ldots
,d_n$ that for all $n$ sufficiently large
\[
\PP\bigl(\bigl\Vert \langle u_{n,i}, f\rangle- \lambda_{n,i}W_{(n)}^{(i)}/c_n
\bigr\Vert < \varepsilon\bigr) \leq\exp \biggl(-I\bigl(\langle u_{n,i}, f
\rangle_{\varepsilon
}\bigr)\frac{
c_n^2}{2n\lambda_{n,i}^2} \biggr), \qquad 1 \le i \le
d_n.
\]
This means that (\ref{51}) and (\ref{53}) are equivalent.
\end{pf}

To further simplify the above criterion for clustering we need the
following uniform lower semicontinuity property of the $I$-function.
%
\begin{lem} \label{lem500}
Let $f=(f_1,\ldots,f_d)$ be such that $\sum_{j=1}^d I(f_j)< \infty$
and $\delta>0$. Then, there exists $\varepsilon>0$ sufficiently small
such that
%
%
\begin{equation}
\label{100} \bigl(I^{{1}/{2}}\bigl(\langle u,f \rangle\bigr) - \delta
\bigr)_{+}^2 \le I\bigl(\langle u,f\rangle_{\varepsilon}
\bigr)
\end{equation}
for all $u \in U=\{u\dvtx |u| \le1\}$.
\end{lem}

\begin{pf}
Let $U_{\delta,f}= \{ u \in U\dvtx I^{{1}/{2}}(\langle u,f
\rangle) \ge\delta\}$. Since $f$ is fixed, $I(\langle u,f\rangle)$ is
continuous and nonnegative on $U$, and the set $U_{\delta,f}$ is
compact. Furthermore, for all $ u \in U \cap U_{\delta,f}^c$ the
conclusion in (\ref{100}) is obvious.
Therefore, if (\ref{100}) fails, it must fail on $ U_{\delta,f}$ and
there exists $u_n \in U_{\delta,f}$ such that for all $n$ sufficiently large
%
%
\begin{equation}
\label{200} \bigl(I^{{1}/{2}}\bigl(\langle u_n,f \rangle\bigr) -
\delta\bigr)_{+}^2 > I\bigl(\langle u_n,f
\rangle_{1/n}\bigr).
\end{equation}
Since $ U_{\delta,f}$ is compact, there is a subsequence $\{u_{n_{k}}\}
$ in $U_{\delta,f}$ and $u_0 \in U_{\delta,f}$ such that $u_{n_{k}}$
converges to $u_0$ and (\ref{200}) holds for $n=n_k$, $k\ge1$.
Using the continuity of $\langle u,f\rangle$ and $I(\langle u,f\rangle
)$ again, we thus have from the left term in (\ref{200}) that
%
%
\begin{equation}
\label{300} \lim_{k \rightarrow\infty} \bigl(I^{{1}/{2}}\bigl(\langle
u_{n_{k}},f \rangle \bigr) - \delta\bigr)_{+}^2 =
\bigl(I^{{1}/{2}}\bigl(\langle u_{0},f \rangle\bigr) - \delta
\bigr)_{+}^2.
\end{equation}
Moreover, since $\langle u,f\rangle$ is continuous on $U$, we have
\[
\lim_{k \rightarrow\infty}\bigl \Vert \langle u_{n_{k}},f
\rangle_{1/n_{k}}- \langle u_0,f\rangle\bigr\Vert =0.
\]
Since the $I$ function is lower semi-continuous and nonnegative, it
follows that
%
%
\begin{equation}
\label{400} \liminf_{k \rightarrow\infty} I\bigl(\langle u_{n_{k}},f
\rangle_{1/n_{k}}\bigr) \ge I\bigl(\langle u_0,f \rangle\bigr).
\end{equation}
Hence, combining (\ref{200}), (\ref{300}) and (\ref{400}) we get
\[
\bigl(I^{{1}/{2}}\bigl(\langle u_{0},f \rangle\bigr) - \delta
\bigr)_{+}^2 \ge I\bigl(\langle u_0,f \rangle
\bigr),
\]
which is a contradiction since $u_0 \in U_{\delta,f}$. Hence the lemma
is proven.
\end{pf}

We can now prove another lemma which will be the crucial tool for
establishing Theorems \ref{theo21} and \ref{theo23}.
%
\begin{lem} \label{lem51}
Under the assumptions of Theorem \ref{theo21} we have
$f=(f_1,\ldots,\break  f_d) \in\mathcal{A}$ if and only if
%
%
\begin{equation}
\label{51a} \sum_{n \in\mathbb{N}_{\varepsilon}} n^{-1} \exp
\Biggl(-\sum_{i=1}^{d_n}\frac{(I^{1/2}(\langle u_{n,i}, f\rangle) -\varepsilon)_{+}^2
c_n^2}{2n\lambda_{n,i}^2}
\Biggr) =\infty \qquad\forall\varepsilon>0,
\end{equation}
where $\mathbb{N}_{\varepsilon}=\{n \ge1\dvtx \|\langle u_{n,i}, f\rangle
\| <
\varepsilon, i >d_n\}$.

Furthermore, we have $x=(x_1,\ldots,x_d) \in A$ if and only if
%
%
\begin{equation}
\label{52} \sum_{n \in\mathbb{N}'_{\varepsilon}} n^{-1} \exp \Biggl(-
\sum_{i=1}^{d_n}\frac{(|\langle u_{n,i}, x\rangle| -{\varepsilon})_{+}^2
c_n^2}{2n\lambda_{n,i}^2} \Biggr) =
\infty\qquad \forall\varepsilon>0,
\end{equation}
where $\mathbb{N}'_{\varepsilon}=\{n \ge1\dvtx |\langle u_{n,i}, x\rangle| <
\varepsilon, i >d_n\}$.
\end{lem}

\begin{pf}
Combining
Lemmas \ref{lem50} and \ref{lem500}, we immediately see that (\ref
{51a}) is necessary for $ f \in\mathcal{A}$.
To show that this condition is also sufficient, it is enough to prove
that (\ref{51a}) implies (\ref{51}); see Lemma \ref{lem50}. To that end
we first note that
since $f=(f_1,\ldots,f_d)$ is fixed and such that $\sum_{j=1}^d I(f_j)<
\infty$, we have $\langle u,f \rangle$ and $I( \langle u,f \rangle)$
both continuous on $U=\{u\dvtx |u|\le1\}$. In addition,
$\langle u,f \rangle_{\varepsilon}$ is jointly continuous in $(\varepsilon
,u)$ with the product topology on $(0,\infty)\times U$ and either the
sup-norm topology or the $H$-norm topology on the range space; see, for
instance, Proposition~2, parts (a) and (b), in \cite{kzi}.

Hence fix $\theta>0$, and set $E_{\theta}=\{u \in U\dvtx \Vert \langle u,f
\rangle\Vert  \ge\theta\}$. We claim that there exists a $\delta>0$
sufficiently small such that
%
%
\begin{equation}
\label{1000} I^{{1}/{2}}\bigl( \langle u,f \rangle_{\theta}\bigr) \le
I^{{1}/{2}}\bigl( \langle u,f \rangle\bigr) - \delta\qquad \forall u \in
E_{\theta}.
\end{equation}
Since $I^{{1}/{2}}( \langle u,f \rangle)$ is continuous on $U$ we
have that $E_{\theta}$ is a compact subset of $U$. Moreover, for $u
\in
E_{\theta}$ we have $I^{{1}/{2}}( \langle u,f \rangle) \ge
\Vert \langle u,f \rangle\Vert  \ge\theta>0$, and consequently
$
I^{{1}/{2}}( \langle u,f \rangle_{\theta}) < I^{{1}/{2}}(
\langle u,f \rangle)
$.

Next define for $k \ge1$,
\[
V_k=\bigl\{u \in U\dvtx I^{{1}/{2}}\bigl( \langle u,f
\rangle_{\theta}\bigr) <I^{{1}/{2}}\bigl( \langle u,f \rangle\bigr) -
1/k \bigr\}.
\]
Then $V_k$ is open by the continuity properties mentioned above, and
$E_{\theta}=\bigcup_{k\ge1}V_k$,\vspace*{2pt} so the compactness of $E_{\theta}$
implies $E_{\theta} \subset V_{k_{0}}$ for some $k_0< \infty$.
Thus (\ref{1000}) holds for $u \in E_{\theta}$ for $\delta=
1/k_{0}$.

If $u \in U\cap E_{\theta}^c$, then we have trivially, $
I^{{1}/{2}}( \langle u,f \rangle_{\theta})=0$. Combining this with
relation (\ref{1000}) and setting $\theta= \varepsilon$, we can conclude
that uniformly on $U$,
\[
I\bigl( \langle u,f \rangle_{\varepsilon}\bigr) \le\bigl(I^{{1}/{2}}\bigl(
\langle u,f \rangle\bigr) - \delta\bigr)_{+}^2,
\]
and we see that indeed (\ref{51a}) implies (\ref{51}).

To prove the second part of Lemma \ref{lem51} we conclude by an obvious
modification of the argument used in Lemma \ref{lem50} that $x \in A$
if and only if
%
%
\begin{equation}
\label{55} \sum_{n=1}^{\infty} n^{-1}
\prod_{i=1}^d \PP\bigl(\bigl|
\lambda_{n,i} W^{(i)}(n)/c_n - \langle
u_{n,i}, x\rangle\bigr| < \varepsilon\bigr)=\infty,\qquad \varepsilon>0,
\end{equation}
where $W^{(i)}(n)\stackrel{d}{=} \sqrt{n}Z$ with $Z$ standard normal.
Consequently we have $x \in A$ if and only if
%
%
\begin{equation}
\label{56} \sum_{n \in\mathbb{N}'_{\varepsilon}} n^{-1} \prod
_{i=1}^{d_n} \PP \bigl(\bigl|\lambda_{n,i}
\sqrt{n} Z - c_n \langle u_{n,i}, x\rangle\bigr| < \varepsilon
c_n\bigr)=\infty,\qquad \varepsilon>0.
\end{equation}
Using a standard argument (see, e.g., part (iii) of the proof
of Proposition 1 in \cite{E-1}) we have that (\ref{56}) holds for all
$\varepsilon>0$ if and only if
\[
\sum_{n \in\mathbb{N}'_{\varepsilon}} n^{-1}\exp \Biggl( -\sum
_{i=1}^{d_n} \frac{(|\langle u_{n,i}, x\rangle| - \varepsilon)_+^2c_n^2}{2n\lambda
_{n,i}^2} \Biggr)=
\infty,\qquad \varepsilon> 0.
\]
Therefore, $x=(x_1,\ldots,x_d) \in A$ if and only if (\ref{52}) holds
for all $\varepsilon>0$.
\end{pf}

We are ready to prove (\ref{25}). Take $x=(x_1,\ldots,x_d) \in A$ and
consider the function $g=(x_1,\ldots,x_d)f$, where $f \in\mathcal{K}$.
Then we have for any vector $u \in\R^d$ and $\varepsilon> 0$,
\[
I\bigl(\langle u, g\rangle\bigr) = I\bigl(\langle u, x \rangle f\bigr)= \langle
u, x\rangle ^2 I(f) \le\langle u, x\rangle^2,
\]
which trivially implies for any $\varepsilon>0$,
$(I^{1/2}(\langle u, g\rangle)-\varepsilon)_{+} \le(|\langle u,x\rangle|
- \varepsilon)_{+}$.

Finally noting that $\mathbb{N}_{\varepsilon} \supset\mathbb
{N}'_{\varepsilon
}$ for this choice of $x$ and $g$ (recall that we have $\|f\|\le1, f
\in\mathcal{K}$), we see that the series for $g$ in (\ref{51a}) must
diverge whenever the series for $x$ in (\ref{52}) diverge. This is of
course the case since we are assuming that $x \in A$.
Thus we have by Lemma \ref{lem51} that $g \in\mathcal{A}$ and relation
(\ref{25}) has been proven.

We next show that $\mathcal{A}$ is star-like and symmetric about zero.
Both properties are direct consequences of Lemma \ref{lem51}.
The symmetry of $\mathcal{ A}$ follows since
\[
I\bigl(\langle f, u \rangle\bigr) = I \bigl(\langle-f, u \rangle\bigr),\qquad f \in
C_d[0,1], u \in\R^d.
\]
To prove that $\mathcal{ A}$ is star-like, we use the simple inequality
$I(\langle\lambda f, u\rangle)= \lambda^2 I(\langle f, u\rangle) \le
I(\langle f, u\rangle)$ which holds for $u \in\R^d$, $0 \le\lambda
\le
1$ and $f=(f_1,\ldots,f_d) \in C_d[0,1]$. It is then obvious
that if $f \in\mathcal{A}$ and consequently the series for $f$ in~(\ref{51a}) diverge, the series for $\lambda f$ must diverge as well,
whence $\lambda f \in\mathcal{A}$.

If $f=(f_1,\ldots,f_d) \in\mathcal{ A}$, then $f(t) \in C(\{
S_{(n)}(t)/c_n\})$ for each fixed $t \in[0,1]$. Now by Theorem 2 in
\cite{E-2}, on a suitable probability space, one can construct a
standard Brownian motion $\tilde{W}(t), t \ge0$ so that with probability
%
%
\begin{equation}
\label{57} \limsup_{n \rightarrow\infty} \|S_{(n)}/c_n -
\Gamma_n\tilde {W}_{(n)}/c_n\| = 0.
\end{equation}
Since $f(t) \in C(\{S_{(n)}(t)/c_n\})$,
we can infer that with probability one
%
%
\begin{equation}
\label{58} \liminf_{n \rightarrow\infty} \bigl|f(t) - \Gamma_n\tilde
{W}_{(n)}(t)/c_n\bigr| =0
\end{equation}
for each $t \in[0,1]$.
Using the scaling property of Brownian motion and (\ref{58}), with
$0<t\leq1$, implies with probability one that
%
%
\begin{equation}
\label{59} \liminf_{n \rightarrow\infty} \bigl|f(t)/\sqrt{t} - \Gamma_n
\tilde {W}_{(n)}(1)/c_n\bigr| =0.
\end{equation}
Thus by (\ref{57}) and (\ref{59}) we have $f(t)/\sqrt{t} \in C(\{
S_{(n)}(1)/c_n\}) =A$ for $0<t \leq1$.

Moreover $A$ is star-like about zero as can be seen directly from Lemma
\ref{lem51} or from the fact that $A=\{f(1)\dvtx f \in\mathcal{A}\}$,
where $\mathcal{A}$ is star-like about zero. Therefore $f(t) \in A$ for
$0\leq t \leq1$. Furthermore,
since $f \in\mathcal{A} \subset C_d[0,1]$, we have that $f$ maps
$[0,1]$ continuously into $A$, and Theorem \ref{theo21} is proven.

\section{\texorpdfstring{Proof of Theorem \protect\ref{theo22}}
{Proof of Theorem 2.2}}\label{sec6} We can assume w.l.o.g.
that $\E(X^{(i)})^2 > 0, 1 \le i \le d$ so that
we have for some $n_0 \ge1$,
\[
\sigma^2_{n,i} = \E\bigl(X^{(i)}
\bigr)^2 I\bigl\{\bigl|X^{(i)}\bigr| \le c_n\bigr\} > 0, \qquad 1
\le i \le d, n \ge n_0.
\]
We then have the following analogue of Lemma \ref{lem51}:
%
\begin{lem} \label{lem61}
Under the assumptions of Theorem \ref{theo22} we have
$f=(f_1,\ldots,\break  f_d) \in\mathcal{A}$ if and only if
%
%
\begin{equation}
\label{61} \sum_{n =n_0}^{\infty}
\frac{1}{n} \exp \Biggl(-\sum_{i=1}^{d}
\frac
{(I^{1/2}(f_i)-{\varepsilon})_{+}^2 c_n^2}{2n\sigma_{n,i}^2} \Biggr) =\infty\qquad \forall\varepsilon>0.
\end{equation}
Furthermore, we have $x=(x_1,\ldots,x_d) \in A$ if and only if
%
%
\begin{equation}
\label{62} \sum_{n =n_0}^{\infty}
\frac{1}{n} \exp \Biggl(-\sum_{i=1}^{d}
\frac
{(|x_i| -{\varepsilon})_{+}^2 c_n^2}{2n\sigma_{n,i}^2} \Biggr) =\infty\qquad \forall\varepsilon>0.
\end{equation}
\end{lem}
The proof is omitted since it is similar to that of Lemma \ref{lem51}.
Simply use Theorem \ref{crimod} instead of Theorem \ref{cri} and part
(c) of Theorem \ref{criA} instead of part (b).

We are ready to prove that
\[
\mathcal{A}=\{x_1\mathcal{K} \times\cdots\times x_d
\mathcal{K} \dvtx x \in A\}.
\]
We first establish the inclusion ``$\supset$.'' Take $x=(x_1,\ldots
,x_d) \in A$ and set $f=(x_1 g_1, \ldots, x_d g_d)$,
where $g_i \in\mathcal{K}, 1 \le i \le d$. Then we obviously have
$I(f_i)= x_i^2 I(g_i) \le x_i^2, 1 \le i \le d$, and we see that
\[
\sum_{n=n_0}^{\infty} \frac{1}{n} \exp
\Biggl(-\sum_{i=1}^{d}\frac{(I^{1/2}(f_i)-{\varepsilon})_{+}^2 c_n^2}{2n\sigma
_{n,i}^2}
\Biggr) \ge\sum_{n =n_0}^{\infty}
\frac{1}{n} \exp \Biggl(-\sum_{i=1}^{d}
\frac{(|x_i| -{\varepsilon})_{+}^2 c_n^2}{2n\sigma
_{n,i}^2} \Biggr),
\]
where the last series is divergent since $x \in A$. In view of Lemma
\ref{lem61} this means that $f \in\mathcal{A}$.

To establish the reverse inclusion ``$\subset$,'' take $f=(f_1,\ldots
,f_d) \in\mathcal{ A}$. From Theorem \ref{theo21} we know that $I(f_i)
< \infty, 1 \le i \le d$. Setting $g_i = f_i/ \sqrt{I(f_i)}$, $1 \le i
\le d$, where $g_i=0$ if $I(f_i)=0$, we have $g_i \in\mathcal{K}, 1
\le i \le d$ and $f= (x_1 g_1, \ldots, x_d g_d)$ if $x_i = \sqrt {I(f_i)}, 1 \le i \le d$, and it is enough to show that $x \in A$. This
is trivial with the above choice for $x$ since for any $\varepsilon>0$,
\[
\sum_{n=n_0}^{\infty} \frac{1}{n}\exp
\Biggl(- \sum_{i=1}^d\frac{(|x_i|-
\varepsilon)_+^2c_n^2}{2n\sigma_{n,i}^2}
\Biggr)=\sum_{n=n_0}^{\infty} \frac
{1}{n}
\exp \Biggl(- \sum_{i=1}^d
\frac{(I^{1/2}(f_i)-{\varepsilon})_{+}^2
c_n^2}{2n\sigma_{n,i}^2} \Biggr),
\]
where the second series is divergent.
Therefore, $x \in A$ by Lemma \ref{lem61}, and Theorem \ref{theo22} has
been proven.
\section{\texorpdfstring{Proof of Theorem \protect\ref{theo23}}
{Proof of Theorem 2.3}}\label{sec7}
W.l.o.g. we can assume that there exists an $n_0 \ge1$ so that all
the matrices $\Gamma_n, n \ge n_0$ have full rank which means that we
have in Lemma \ref{lem51} $d_n =2, n \ge n_0$. Otherwise, $X$ is
supported by a 1-dimensional subspace of $\R^2$ (see the comments
before Lemma \ref{lem50}) and in this case it easily follows from the
1-dimensional functional LIL type result in \cite{E} that $\mathcal
{A}=\{xg\dvtx  x \in A, g \in\mathcal{K}\}$ which trivially implies the
assertion of Theorem \ref{theo23}.

We show that any function $f \in\mathcal{A}$ has a representation
$(x_1 g_1, x_2 g_2)$, where $(x_1, x_2) \in A$ and $g_1, g_2 \in
\mathcal{K}$.
To that end we look first at ``nonextremal'' functions $f \in\mathcal
{A}$. That is, we assume that $f \in\mathcal{A}$ is such that
$(1+\eta
)f \in\mathcal A$ for some $\eta>0$. Also assume that $f \ne0$.

Rewrite $f$ as $(x_1 h_1,x_2 h_2)$, where $I(h_i)=1, i =1, 2$.

Note that $I(\langle u, f\rangle)= \sum_{i,j=1}^2 x_i u_i x_j u_j
\alpha
_{i,j}$, where
$\alpha_{i,j} = \int_0^1 h'_i(s) h'_j(s) \,ds$ for $ 1 \le i,j \le2$.
Then we obviously have $\alpha_{1,1}=\alpha_{2,2}=1$ and consequently
\[
I\bigl(\langle u, f\rangle\bigr) = x_1^2
u_1^2 +x_2^2 u_2^2
+ 2\alpha_{1,2} u_1x_1 u_2
x_2,
\]
where $|\alpha_{1,2}| \le1$ by Cauchy--Schwarz.

Similarly, we have for $y \in\R^2$,
\[
\langle u, y\rangle^2 = y_1^2
u_1^2 +y_2^2 u_2^2
+ 2 u_1y_1 u_2 y_2.
\]
Set $\tilde{x}=(-x_1, x_2)$. Comparing the two expressions above we
see that
%
%
\begin{equation}
\label{90a} I\bigl(\langle u, f\rangle\bigr) \ge\langle u, x\rangle^2
\wedge\langle u, \tilde{x}\rangle^2.
\end{equation}
Next observe that we have if $\{u_{n,1}, u_{n,2}\}$ is an orthonormal
basis of $\R^2$,
%
%
\begin{equation}
\label{91} \sum_{i=1}^2 I\bigl(\langle
u_{n,i}, f\rangle\bigr)= \int_0^1
\sum_{i=1}^2 \bigl\langle u_{n,i},
f'(s)\bigr\rangle^2 \,ds= \int_0^1
\bigl|f'(s)\bigr|^2\,ds=|x|^2.
\end{equation}
Further note that $(1+\eta)f \in\mathcal{A}$ implies via Lemma \ref
{lem51} that
%
%
\begin{equation}\qquad
\sum_{n=n_0}^{\infty} n^{-1} \exp
\Biggl(-\sum_{i=1}^{2}\frac
{(I^{1/2}(\langle u_{n,i}, (1+\eta) f\rangle)-\varepsilon)_{+}^2
c_n^2}{2n\lambda_{n,i}^2}
\Biggr) =\infty,\qquad\varepsilon>0.
\end{equation}
In view of relation (\ref{91}) we can find a sequence $i_n \in\{1,2\}$
so that
\[
I\bigl(\langle u_{n,i_n}, f\rangle\bigr) \ge|x|^2/2,\qquad n \ge1.
\]
Then one must have
\[
\sum_{n\dvtx i_n=1} n^{-1} \exp \Biggl(-\sum
_{i=1}^{2}\frac{((1+\eta
)I^{1/2}(\langle u_{n,i}, f\rangle)-\varepsilon)_{+}^2 c_n^2}{2n\lambda
_{n,i}^2} \Biggr) =
\infty,\qquad \varepsilon>0
\]
or
\[
\sum_{n: i_n=2} n^{-1} \exp \Biggl(-\sum
_{i=1}^{2}\frac{((1+\eta
)I^{1/2}(\langle u_{n,i}, f\rangle)-\varepsilon)_{+}^2 c_n^2}{2n\lambda
_{n,i}^2} \Biggr) =
\infty, \qquad\varepsilon>0.
\]
We can assume w.l.o.g. that the series for $i_n=1$ diverge. Then an
easy calculation shows that if $0 < \varepsilon< \eta|x|/\sqrt{2}$,
\[
\sum_{n=n_0}^{\infty} n^{-1} \exp
\biggl(-\frac{I(\langle u_{n,1},
f\rangle) c_n^2}{2n\lambda_{n,1}^2} - \frac{((1+\eta
)I^{1/2}(\langle
u_{n,2}, f\rangle)-\varepsilon)_{+}^2 c_n^2}{2n\lambda_{n,2}^2} \biggr) =\infty.
\]
Next set for $\beta>0$,
\[
J(\beta)=\bigl\{n \ge n_0\dvtx I\bigl(\langle u_{n,2}, f
\rangle\bigr) \le\beta^2\bigr\}
\]
and
\[
\rho= \inf \biggl\{ \beta> 0\dvtx \sum_{n \in J(\beta)}
n^{-1} \exp \biggl(-\frac{I(\langle u_{n,1}, f\rangle) c_n^2}{2n\lambda_{n,1}^2} \biggr)=\infty \biggr\}.
\]
There are two cases:

\textit{Case} 1 \fbox{$\rho> 0$}. We then can choose an arbitrary $0 <
\rho
_1 < \rho$, and we get that
\[
\sum_{n \notin J(\rho_1)} n^{-1} \exp \biggl(-
\frac{I(\langle u_{n,1},
f\rangle) c_n^2}{2n\lambda_{n,1}^2} - \frac{((1+\eta
)I^{1/2}(\langle
u_{n,2}, f\rangle)-\varepsilon)_{+}^2 c_n^2}{2n\lambda_{n,2}^2} \biggr) =\infty.
\]
Noticing that $I^{1/2}(\langle u_{n,2}, f\rangle) \ge\rho_1$ if $n
\notin J(\rho_1)$, we can conclude if $\varepsilon< \rho_1 \eta$ that
\[
\sum_{n \notin J(\rho_1)} n^{-1} \exp \Biggl(-\sum
_{i=1}^2 \frac
{I(\langle u_{n,i}, f\rangle) c_n^2}{2n\lambda_{n,i}^2} \Biggr) =
\infty.
\]
Set $\mu_{n,1} = \lambda_{n,1} \vee\lambda_{n,2}$, $\mu_{n,2} =
\lambda
_{n,1} \wedge\lambda_{n,2}$, and denote the corresponding eigenvectors
in $\{u_{n,1}, u_{n,2}\}$ by $v_{n,1}$ and $v_{n,2}$.
Then we have by (\ref{91}),
\[
\sum_{i=1}^{2} I\bigl(\langle
u_{n,i}, f\rangle\bigr)/\lambda_{n,i}^2
=|x|^2/\mu_{n,1}^2 + \bigl(
\mu_{n,2}^{-2}-\mu_{n,1}^{-2}\bigr)I\bigl(
\langle v_{n,2}, f\rangle\bigr),
\]
where $\mu_{n,2}^{-2}-\mu_{n,1}^{-2} \ge0$.

In view of (\ref{90a}) we can find a sequence $a_n \in\{-1,1\}$ so
that we have for $y_n=(a_n x_1, x_2)$,
\[
\langle v_{n,2}, y_n\rangle^2 \le I\bigl(
\langle v_{n,2}, f\rangle\bigr),\qquad n \ge1,
\]
which then implies that
\[
\sum_{i=1}^{2}I\bigl(\langle
u_{n,i}, f\rangle\bigr)/\lambda_{n,i}^2
\ge|y_n|^2/\mu_{n,1}^2 + \bigl(
\mu_{n,2}^{-2}-\mu_{n,1}^{-2}\bigr)\langle
v_{n,2}, y_n\rangle^2 =\sum
_{i=1}^{2}\langle u_{n,i}, y_n
\rangle^2/\lambda_{n,i}^2.
\]
It follows that
\[
\sum_{n=n_0}^{\infty} n^{-1} \exp
\Biggl(-\sum_{i=1}^{2}\frac
{\langle
u_{n,i}, y_n\rangle^2 c_n^2}{2n\lambda_{n,i}^2}
\Biggr) =\infty.
\]
But this implies that
\begin{eqnarray*}
\sum_{n\dvtx a_n=1} n^{-1} \exp \Biggl(-\sum
_{i=1}^{2}\frac{\langle u_{n,i},
x\rangle^2 c_n^2}{2n\lambda_{n,i}^2} \Biggr)&=&\infty\quad
\mbox{or}
\\
\sum_{n\dvtx a_n=-1} n^{-1} \exp \Biggl(-\sum
_{i=1}^{2}\frac{\langle
u_{n,i}, \tilde{x}\rangle^2 c_n^2}{2n\lambda_{n,i}^2} \Biggr)&=&\infty.
\end{eqnarray*}
Recalling Lemma \ref{lem51} we see that $f=(x_1 h_1, x_2 h_2) \in
\mathcal{A}$ implies $(x_1, x_2) \in A$ or $\tilde{x}=(-x_1, x_2) \in
A. $

Rewriting $f$ as $(-x_1 g_1, x_2 g_2)$ if $\tilde{x} \in A$,
where $g_1=-h_1, g_2=h_2$ we see that $f$ has always the desired form
in Case 1.

\textit{Case} 2 \fbox{$\rho=0$}. In this case we have by definition of
$\rho
$ for any $\varepsilon>0$,
\[
\sum_{n \in J(\varepsilon)} n^{-1} \exp \biggl(-
\frac{I(\langle u_{n,1},
f\rangle) c_n^2}{2n\lambda_{n,1}^2} \biggr)=\infty,
\]
which in turn implies if $\varepsilon< |x|$ that
\[
\sum_{n \in J(\varepsilon)} n^{-1} \exp \biggl(-
\frac{(|x|^2 -\varepsilon^2)
c_n^2}{2n\lambda_{n,1}^2} \biggr)=\infty.
\]
Here we have again used relation (\ref{91}) from which we can infer
that
\[
I\bigl(\langle u_{n,1}, f\rangle\bigr) \ge|x|^2 -
\varepsilon^2,\qquad n \in J(\varepsilon).
\]
Choosing $y_n=(\pm x_1, x_2)$ so that $\langle u_{n,2}, y_n \rangle^2
\le I(\langle u_{n,2}, f \rangle), n \ge1$,
we get for $\varepsilon< |x|$ and $n \in J(\varepsilon)$,
\[
\sum_{i=1}^2 \frac{(|\langle u_{n,i}, y_n \rangle| - \varepsilon
)_{+}^2}{\lambda_{n,i}^2} =
\frac{(|\langle u_{n,1}, y_n \rangle| -
\varepsilon)_{+}^2}{\lambda_{n,1}^2} \le\frac{(|y_n| -\varepsilon
)^2}{\lambda
_{n,1}^2} \le\frac{|x|^2 -\varepsilon^2}{\lambda_{n,1}^2},
\]
and we can conclude that
\[
\sum_{n = n_0}^{\infty} n^{-1} \exp
\Biggl(-\sum_{i=1}^2 \frac
{(|\langle
u_{n,i}, y_n \rangle| - \varepsilon)_{+}^2 c_n^2}{2n\lambda
_{n,i}^2}
\Biggr)=\infty,\qquad  \varepsilon< |x|.
\]
This implies as in Case 1 that $(x_1,x_2) \in A$ or $\tilde
{x}=(-x_1,x_2) \in A$ and finally that $f$ has the desired form.

If $f$ is an extremal function we can find a sequence $f_n$ of
nonextremal functions converging to it (in sup-norm). These functions
$f_n$ have the form $(x_{n,1} g_{n,1}, x_{n,2} g_{n,2})$ where
$(x_{n,1}, x_{n,2})\in A$
and $ g_{n,i} \in\mathcal{K}, i=1,2$. By compactness of $A$ and
$\mathcal{K}$ we can find a subsequence $n_k$
so that $(x_{n_k,1}, x_{n_k,2})$ and $g_{n_k,i}$ converge to $(x_1,x_2)
\in A$ and $g_i \in\mathcal{K}, i=1,2$, respectively. Consequently we
have $f=\lim_{k\to\infty} (x_{n_k,1} g_{n_k,1}, x_{n_k,2}
g_{n_k,2})=(x_1g_1, x_2 g_2)$ and Theorem \ref{theo23} has been proven.

\begin{rems*}
(1) The same proof shows that if we
use an arbitrary orthonormal basis $\{u, v\}$ of $\R^2$ to express $X$,
then we have 
\[
\mathcal{A} \subset\bigl\{f_1\langle x, u \rangle u +
f_2\langle x, v \rangle v\dvtx f_1, f_2 \in
\mathcal{K}, x \in A\bigr\}.
\]
In certain cases this can lead to a smaller upper bound set than that
one obtained from Theorem \ref{theo23}, which has $X$ given in terms of
the canonical basis.

(2) One might wonder whether the result also holds in dimension $d \ge
3$. In the present proof we have used the following fact about
quadratic forms in $\R^2$ [see~(\ref{90a})] which has no direct
analogue in higher dimensions:
Given two symmetric positive semidefinite $(2,2)$-matrices $A, B$ with
$A_{i,i} = B_{i,i}, i=1,2$ and $|A_{1,2}| \le|B_{1,2}|$, one has for
any $x=(x_1,x_2) \in\R^2$: $\langle x, A x\rangle\ge\langle x, B
x\rangle\wedge\langle\tilde{x}, B \tilde{x}\rangle$,
where $\tilde{x}=(-x_1, x_2)$.
\end{rems*}

So clearly a different proof would be necessary in order to prove this
result in higher dimensions if this is possible at all.

\section{An example}\label{sec8}
In this final section we show that for any nonempty closed subset
$\tilde{A}$ of $\
R^d$ which is star-like and symmetric w.r.t. 0 there are
$d$-dimensional distributions such that $\tilde{A}$ is the cluster set
for $S_n/c_n$ and
at the same time the functional cluster set $\mathcal{A}$ is of the
form $\{xg\dvtx  x\in\tilde{A}, g \in\mathcal{K}\}$.

This can be
done for the generalized LIL in \cite{EL-1}; that is, such
distributions exist for the normalizing sequence $c_n =\sqrt{2n (\log
\log n)^{1+p}}$, where $p>0$. To simplify notation, we will prove this
only if $p=1$ and if the set $\tilde{A}$ is such that $\max_{x \in
\tilde{A}} |x| =1$. It should be obvious to the reader how to do the
``general'' case once he or she has seen how it works for this special
case.

The point is that this phenomenon occurs for very regular
normalizing sequences.
%
\begin{theo}\label{th8.1}
Let $\tilde{A}$ be a set in $\R^d$ which is symmetric and star-like
with respect to zero and which satisfies $\max_{x \in\tilde{A}} |x| =1$.
Then, one can find a $d$-dimensional distribution $Q$ such that
for $X_1, X_2, \ldots$ independent
$Q$-distributed random vectors
and $S_n=\sum_{j=1}^n X_j, n \ge1$, we have with probability one,
%
%
\begin{eqnarray}
\limsup_{n \to\infty} |S_n|/\sqrt{2n}(\log\log n)& =&1,
\label {ex1}
\\
C\bigl(\bigl\{S_n/\sqrt{2n} (\log\log n)\dvtx n \ge3\bigr\}\bigr)&=&
\tilde{A},\label{ex2}
\\
C\bigl(\bigl\{S_{(n)}/\sqrt{2n} (\log\log n)\dvtx n \ge3\bigr\}
\bigr)&=&\bigl\{(x_1 g, \ldots, x_d g) \dvtx g \in
\mathcal{K}, x \in\tilde{A} \bigr\}.\label{ex3}
\end{eqnarray}
\end{theo}
To prove this result, we use a similar idea as in Theorem 5 of \cite
{E-1} and Theorem 2 of \cite{EK}: we start with the construction of a
real random variable $Z$ in the domain of attraction of the normal
distribution, and then we define a suitable random vector $X\dvtx \Omega
\to\R^d$ as a function of this variable $Z$. Due to the use of the
normalizing sequence $c_n =\sqrt{2n}\log\log n$ instead of the
normalizers used in \cite{E-1,EK}, and the recent work of \cite{EL-1,EL}, some simplification is possible.

\begin{pf*}{Proof of Theorem \ref{th8.1}}
\textit{Step} 1. \textit{Definition of the random variable $Z$.}
We first define a monotone right continuous function $H\dvtx [0,\infty[
\,\to
[0,\infty[$ which satisfies
\[
\liminf_{t \to\infty} H(t)/\log\log t = 0 \quad\mbox{and}\quad \limsup
_{t
\to
\infty} H(t)/\log\log t = 1.
\]
We set for $k \geq1$, $m_k = 3^{2^{k^3}}$, $m_{k,0} = m_k$ and
$m_{k,\ell
}= 3^{2^{k^3 + \ell k}}$ for $0 \leq\ell
\leq k$. Furthermore, we define $m_{k,k+1} = m_{k+1}$ and
$n_{k,\ell} = m_{k,\ell+1} - k^3, 0 \leq\ell\leq k$.

We assume that $H(t), t \geq0$ satisfies
\[
H(t) = d_n,\qquad \exp( n) \leq t < \exp( n+1), \qquad n \geq1,
\]
where
\begin{eqnarray*}
d_n & =& 0,\qquad 0 \leq n < m_1
\quad\mbox{and for } k \geq1,
\\
d_n &=& d_{m_{k,\ell}} = (\log3)2^{k^3 + \ell k},\qquad
m_{k,l} \le n \le n_{k,l}, 0 \leq\ell \leq k,
\\
d_{n_{k,\ell} +j} & =& (\log3)2^{k^3 + \ell k + j/k^2},\qquad 1 \le j \le k^3, 0
\le\ell\le k-1,
\\
d_{n_{k,k} + j} & =& (\log3) 2^{(2k^2 + 3k + 1) jk^{-3} + k^2 + k^3},\qquad 1 \leq j \leq k^3.
\end{eqnarray*}
%
From this definition we note $d_n$ is defined for every integer $n\ge
0$, and we readily obtain $H(t) \le\log\log t, t \ge e$. We also have
\[
H\bigl(\exp(m_{k,\ell})\bigr) = \log m_{k,\ell},\qquad 0 \le\ell\le k+1,
k \ge1,
\]
so that indeed $\limsup_{t \to\infty} H(t)/\log\log t = 1$.\

Further note that
\[
H(t) = H\bigl(\exp(m_{k,\ell})\bigr),\qquad \exp(m_{k,\ell}) \le t <
\exp (n_{k,\ell
}+1),\qquad 0 \le\ell\le k,
\]
which implies $\liminf_{t \to\infty} H(t)/\log\log t = 0$ as $\log
n_{k,\ell}/ \log m_{k,\ell} \ge2^k-1$ for $0 \le\ell\le k$.

Similarly as in Lemma 8 of \cite{E-1} we define a symmetric and
discrete random variable $Z\dvtx \Omega\to\R$ with support $\{0, \pm
\exp
(n)\dvtx n \ge m_1\}$ such that $\E[Z^2I\{|Z| \le t\}] = H(t), t \ge0$.

To accomplish this we set $q_n = (d_n - d_{n-1}) e^{-2n}/2, n \ge m_1$.

It is easily checked that $\sum_{n=m_1}^{\infty} q_n < 1/2$. Thus
there exists a discrete random variable satisfying
$\PP\{Z= \exp(n)\} = \PP\{Z= -\exp(n)\} = q_n, n \ge m_1$ and $\PP
\{Z=
0\} = 1 -2\sum_{n=m_1}^{\infty} q_n$.

An easy calculation then shows
that $\E[Z^2I\{|Z| \le t\}] = H(t), t \ge0$.

Moreover, since
$d_{n+1}/d_n \to1$ as $n \to\infty$, we have $H(et)/H(t) \to1$ as $t
\to\infty$. Consequently, the function $H$ is slowly varying at
infinity. It follows that $Z$ is in the
domain of attraction of the normal distribution. Recall that this
implies among other things that
%
%
\begin{equation}
\label{exa1} t^2 \PP\bigl\{|Z| > t\bigr\}/H(t) \to0\qquad \mbox{as }t \to\infty.
\end{equation}

\textit{Step} 2. \textit{Definition of the random vector $X\dvtx \Omega\to\R^d.$}
We write the set $\tilde{A}$ as a closure of a union of countably many
symmetric line segments, that is,
$\tilde{A} = \operatorname{cl} (\bigcup_{j=1}^{\infty} \mathcal
{L}_j )$,
where $\mathcal{L}_j=\{tz_j\dvtx |t| \le\sigma_j\}, |z_j|=1$ and $0 <
\sigma_j \le1, j \ge1$.
Note that we also have this representation if $\tilde{A}$ is a union of
finitely many symmetric line segments $\mathcal{L}_j, 1 \le j \le m$.
In this case we simply set $\mathcal{L}_j=\mathcal{L}_1, j \ge m+1$.
Moreover, by repeating some of the line segments $\mathcal{L}_j$ in the
representation of $\tilde A$ if necessary, we may assume without loss
of generality that $\sigma_j^2 \ge1/j, j \ge1$.
Furthermore, we can and do assume that $\sigma_1 =1$ (since there must
be a line segment with $\sigma_j =1$ as $\max_{x \in\tilde{A}} |x| =
1)$.

We now can define a suitable random vector $X\dvtx \Omega\to\R^d$ as follows:
\[
X = \sum_{k=1}^{\infty} \sum
_{\ell=1}^{k+1} \sigma_{\ell} z_{\ell
}
Z I\bigl\{\exp(m_{k,\ell-1}) < |Z| \le\exp(m_{k,\ell})\bigr\}.
\]
From the definition of $X$ and the $H$-function it follows that for any
$x \in\R^d$ with $|x|=1$ and for $\exp(m_{k,\ell-1}) \le t \le\exp
(m_{k,\ell}), 1 \le\ell\le k+1, k \ge1$
%
%
\begin{eqnarray}
\label{exa2} \E\bigl[\langle X, x \rangle^2 I\bigl\{|Z| \le t\bigr\}\bigr] &
\le& H\bigl(\exp(m_{k,\ell
-1})\bigr)
\\
&&{}+ \sigma_{\ell}^2 \langle x, z_{\ell}
\rangle^2 \bigl[H(t)-H\bigl(\exp(m_{k,\ell-1})\bigr)\bigr],
\nonumber
\\
\E\bigl[\langle X, x \rangle^2 I\bigl\{|Z| \le t
\bigr\}\bigr] &\ge& \sigma_{\ell}^2 \langle x, z_{\ell}
\rangle^2 \bigl[H(t)-H\bigl(\exp(m_{k,\ell-1})\bigr)\bigr].\label{e8.6}
\end{eqnarray}
\textit{Step} 3. \textit{Proof of the upper bound in} (\ref{ex1}).

Note
that $X$ has a symmetric distribution since the distribution of $Z$ is
symmetric.
Moreover, $|X| \le|Z|$ so that $\E X $ exists and it has to be equal
to zero by symmetry.

Set
\[
H_X(t):=\sup\bigl\{\E\bigl[\langle v, X\rangle^2 I
\bigl\{|X| \le t\bigr\}\bigr] \dvtx |v| \le1\bigr\},\qquad t \ge0.
\]
Observe that we have for any vector $v \in\R^d$ with $|v| \le1$,
\begin{eqnarray*}
\E\bigl[\langle v, X\rangle^2 I\bigl\{| X| \le t\bigr\}\bigr] &=& \E\bigl[
\langle v, X\rangle^2 I\bigl\{ |Z| \le t\bigr\}\bigr]
\\
&&{} +\E\bigl[\langle v, X\rangle^2 I\bigl\{| X| \le t, |Z| >t\bigr\}\bigr]
\\
&\le&\E\bigl[Z^2I\bigl\{|Z| \le t\bigr\}\bigr] + t^2 \PP\bigl\{|Z| > t
\bigr\}
\\
&=& H(t) + t^2 \PP \bigl\{ |Z| > t\bigr\}.
\end{eqnarray*}
Recalling (\ref{exa1}), we can conclude that $\limsup_{t\to\infty}
H_X(t)/H(t) \le1$ which in turn implies that
\[
\limsup_{t \to\infty} H_X(t)/\log\log t \le1,
\]
and hence by using Corollary 2.5 in \cite{E-2} with $p=2$, we see that
with probability one,
%
%
\begin{equation}
\label{ex0} \limsup_{n \to\infty} |S_n|/\sqrt{2n}\log\log n
\le1.
\end{equation}
%
Moreover, the comment following Corollary 2.5 in \cite{E-2} also implies
%
%
\begin{equation}
\label{exxx} \sum_{n=1}^{\infty}\mathbb{P}\bigl(|X|
\ge\sqrt{2n} {\log\log}n\bigr)<\infty.
\end{equation}
%
\textit{Step} 4. \textit{Proof of} (\ref{ex2}).
Here it is enough to show ``$\supset$.'' The inclusion ``$\subset$''
will follow from the inclusion ``$\subset$'' in (\ref{ex3}) using again
the fact that $A=C(\{S_n/\sqrt{2n} \log\log n\dvtx n \ge3\})$ is equal to
$\{f(1)\dvtx f \in\mathcal{A}\}$ where $\mathcal{A}$ is the corresponding
functional cluster set. To simplify notation, we write $c_n$ instead of
$\sqrt{2n} \log\log n, n \ge3$.

Since $C(\{S_n/c_n\dvtx n \ge3\})$ is
a closed subset of $\R^d$, it is obviously enough to show that
\[
\mathcal{L}_j \subset C\bigl(\{S_n/c_n\dvtx
n \ge3\}\bigr),\qquad j \ge1.
\]
We also know that the cluster set is symmetric and star-like with
respect to 0 so that we only need to prove that
%
%
\begin{equation}
\label{ex4} \sigma_j z_j \in C\bigl(
\{S_n/c_n\dvtx n \ge3\}\bigr), \qquad j \ge1.
\end{equation}
Furthermore, it follows by a slight modification of the proof of
Theorem \ref{criA}(b) [using a different truncation level which is
possible since $|X| \le|Z|$ and $\sum_{n=1}^{\infty} \PP\{|Z| \ge
c_n\} < \infty$
by (\ref{exxx})],
that
%
%
\begin{eqnarray}
\label{ex4a} x &\in& C\bigl(\{S_n/c_n\}\bigr)\qquad \mbox{ a.s.}\quad
\Longleftrightarrow
\nonumber
\\[-8pt]
\\[-8pt]
\nonumber
\sum_{n=3}^{\infty
}
n^{-1} \PP\bigl\{ |x - \sqrt{n} Y_n/c_n| <
\varepsilon\bigr\} &=& \infty,\qquad \varepsilon>0,
\end{eqnarray}
where $Y_n$ is $\operatorname{normal}(0,(\Gamma'_n)^2)$-distributed with $(\Gamma
'_n)^2 = \operatorname{cov}(XI\{|Z| \le c_n\})$.

Set for $k \ge j$,
\[
I_{k,j} =\bigl\{n\dvtx \exp(m_{k,j}) \le c_n \le
\exp(2m_{k,j})\bigr\}.
\]
Then we have by (\ref{e8.6})
\begin{eqnarray*}
\E\bigl[\langle X, z_j \rangle^2 I\bigl\{|Z| \le
c_n\bigr\}\bigr] &\ge& \E\bigl[\langle X, z_j
\rangle^2 I\bigl\{|Z| \le\exp(m_{k,j})\bigr\}\bigr]
\\
&\ge&\sigma_j^2 \bigl[H\bigl(\exp(m_{k,j})
\bigr)- H\bigl(\exp(m_{k,j-1})\bigr)\bigr].
\end{eqnarray*}
Since $H(\exp(m_{k,\ell}))=\log m_{k,\ell}, \ell=j-1, j$ and $\log
m_{k,j} \ge2^k \log m_{k,j-1}$, we get for
$n \in I_{k,j}$,
%
%
\begin{equation}
\label{ex5} \operatorname{Var}\bigl(\langle z_j, Y_n\rangle
\bigr) = \E\bigl[\langle X, z_j \rangle^2 I\bigl\{ |Z| \le
c_n\bigr\}\bigr] \ge\sigma_j^2\bigl(1 -
2^{-k}\bigr)\log m_{k,j}.
\end{equation}
Similarly, we can infer from (\ref{exa2}) for any vector $w$ such that
$|w|=1$ and $\langle w, z_j \rangle=0$,
%
%
\begin{equation}
\label{ex6} \operatorname{Var}\bigl(\langle w, Y_n\rangle\bigr) \le H\bigl(
\exp(m_{k,j-1})\bigr) \le 2^{-k}\log m_{k,j}, \qquad n \in
I_{k,j}.
\end{equation}
Let $0 < \varepsilon< \sigma_j$ and recall that $\sigma_j \le1$.
Choosing an orthonormal basis $\{w_{j,1},\ldots,\break  w_{j,d}\}$ of $\R^d$
with $w_{j,1}=z_j$, we then have with $\varepsilon_1:=\varepsilon/\sqrt{d}$,
\begin{eqnarray*}
&&\PP\{ |\sigma_j z_j - \sqrt{n} Y_n/c_n|
< \varepsilon\}
\\
&&\qquad\ge \PP\bigl\{ \bigl|\sigma_j - \sqrt{n} \langle z_j,
Y_n\rangle/c_n\bigr| < \varepsilon_1, \sqrt{n}\bigl |
\langle w_{j,i}, Y_n\rangle/c_n\bigr| < \varepsilon
_1, 2 \le i \le d \bigr\}
\\
&&\qquad\ge\PP\bigl\{\bigl |\sigma_j - \sqrt{n} \langle z_j,
Y_n\rangle/c_n\bigr| < \varepsilon _1\bigr\} -\sum
_{i=2}^d \PP\bigl\{\bigl|\langle
w_{j,i}, Y_n\rangle\bigr| \ge\varepsilon_1
c_n/\sqrt{n}\bigr\}.
\end{eqnarray*}
Using the fact that $\log\log n \ge\log\log c_n \ge\log m_{k,j}, n
\in I_{k,j}$, we obtain from (\ref{ex6}) for $2 \le i \le d$, $n \in
I_{k,j}$ and large enough $k$,
%
%
\begin{eqnarray}
\label{ex7} \PP\bigl\{\bigl|\langle w_{j,i}, Y_n\rangle\bigr| \ge
\varepsilon_1 c_n/\sqrt{n}\bigr\} &\le& \exp\bigl(-
\varepsilon_1^2 2^k (\log\log
n)^2/\log m_{k,j}\bigr)
\nonumber
\\[-8pt]
\\[-8pt]
\nonumber
&\le& m_{k,j}^{-2}.
\nonumber
\end{eqnarray}
On the other hand, we have for $n \in I_{j,k}$ and large enough $k$,
\[
\PP\bigl\{ \bigl|\sigma_j - \sqrt{n} \langle z_j,
Y_n\rangle/c_n\bigr| < \varepsilon _1\bigr\} \ge
\PP\bigl\{ \langle z_j, Y_n\rangle\ge(
\sigma_j - \varepsilon_1) c_n/\sqrt {n}\bigr
\} /2.
\]
Next observe that
\[
c_n/\sqrt{n}=\sqrt{2}\log\log n \le\sqrt{2} \log\log
c_n^2 \le \sqrt {2}\log(4m_{k,j}), \qquad n \in
I_{k,j}.
\]
Combining (\ref{ex5}) with the obvious fact that $\log(4m_{k,j})/\log
(m_{k,j}) \to1$ as $k \to\infty$, we get for $n \in I_{j,k}$ and
large enough $k$ that
\[
\PP\bigl\{ \langle z_j, Y_n\rangle\ge(
\sigma_j - \varepsilon_1) c_n/\sqrt {n}\bigr
\} \ge\PP\bigl\{\xi> (1-\varepsilon_1 /2)\sqrt{2\log
m_{k,j}}\bigr\},
\]
where $\xi$ is standard normal.

Employing the trivial inequality $\PP\{\xi\ge t\} \ge t^{-1} \exp
(-t^2/2)/\sqrt{8\pi}, t \ge1$,
we get for $n \in I_{k,j}$ and large $k$,
\[
\PP\bigl\{ \bigl|\sigma_j - \sqrt{n} \langle z_j,
Y_n\rangle/c_n\bigr| < \varepsilon _1\bigr\}
\ge(64\pi)^{-1/2}(\log m_{k,j})^{-1}
m_{k,j}^{-1+\varepsilon_1/2}.
\]
Recalling (\ref{ex7}), we can conclude that for $n \in I_{k,j}$ and
large $k$,
\[
\PP\bigl\{ |\sigma_j z_j - \sqrt{n} Y_n/c_n|
< \varepsilon\bigr\} \ge16^{-1} (\log m_{k,j})^{-1}
m_{k,j}^{-1+\varepsilon_1/2}.
\]
Set $a_{k,j}=\min I_{k,j}$ and $b_{k,j}=\max I_{j,k}$. Then we have
\[
\sum_{n \in I_{k,j}}n^{-1} =\sum
_{n = a_{k,j}}^{b_{k,j}}n^{-1} \ge \log
\bigl((b_{k,j}+1)/(a_{k,j}-1)\bigr) -1.
\]
As we have $c_n/c_m \le n/m, n \ge m$, we can infer from the definition
of $I_{k,j}$ that
\[
\sum_{n \in I_{k,j}} n^{-1} \ge m_{k,j}
-1.
\]
We now see that as $k \to\infty$,
\[
\sum_{n \in I_{k,j}}n^{-1} \PP\bigl\{ |
\sigma_j z_j - \sqrt{n} Y_n/c_n|
< \varepsilon\bigr\} \to\infty,
\]
which in view of (\ref{ex4a}) means that $\sigma_j z_j \in C(\{
S_n/c_n\}
), j \ge1$.

Thus $C(\{S_n/\sqrt{2n}\log\log n\dvtx n \ge3\}) \supset\tilde{A}$.

Notice that this also implies
\[
\limsup_{n \to\infty} |S_n|/\sqrt{2n}\log\log n \ge1
\qquad\mbox{a.s.}
\]
since $\tilde{A}$ contains the vector $z_1$ which has norm 1.
Combining this observation with the upper bound (\ref{ex0}) we see
that (\ref{ex1}) holds.

\textit{Step} 5.
\textit{Proof of} (\ref{ex3}).
It remains only to prove the inclusion ``$\subset$.'' The other
inclusion follows directly from the inclusion ``$\supset$'' in (\ref
{ex2}) and (\ref{25}), since the assumptions of Theorem \ref{theo21} hold by
(\ref{exxx}),
and we have already verified (\ref{ex1}).

As in the proof of (\ref{ex2}) we set $c_n=\sqrt{2n}\log\log n, n
\ge
3$. Moreover, we replace the matrices $\Gamma_n$ in Theorem \ref{cri}
by the symmetric and positive semidefinite matrices $\Gamma'_n$ satisfying
\[
\bigl(\Gamma'_n\bigr)^2 = \operatorname{cov}
\bigl(XI\bigl\{|Z| \le c_n\bigr\}\bigr), \qquad n \ge1.
\]
That this is possible follows easily from the proof of Theorem \ref
{cri}. Recall that $|X| \le|Z|$ and that $\sum_{n=1}^\infty\PP\{|Z|
\ge c_n\} < \infty$.

Using again the notation $s_n= S_{(n)}/\sqrt {2n}\log\log n, n \ge3$, we thus have
%
%
\begin{eqnarray}
\label{ex10} f &\in& C\bigl(\{s_n\}\bigr)\qquad \mbox{a.s.}\quad
\Longleftrightarrow
\nonumber
\\[-8pt]
\\[-8pt]
\nonumber
 \sum_{n=3}^{\infty}
n^{-1}\PP\bigl\{\bigl\|\Gamma'_n
W_{(n)}/c_n - f\bigr\| < \varepsilon \bigr\}& =& \infty\qquad\forall
\varepsilon>0.
\end{eqnarray}
Setting $\mathcal{K}_j:= \{\sigma_jz_j f\dvtx f \in\mathcal{K}\}, j
\ge
1$, it is easy to see that
\[
\tilde{\mathcal{A}}:=\{xg\dvtx x \in\tilde{A}, g \in\mathcal{K}\} = \operatorname{cl}
\Biggl(\bigcup_{j=1}^{\infty}
\mathcal{K}_j \Biggr).
\]
We shall show that if $f \notin\mathcal{\tilde{A}}$, then the series
in (\ref{ex10}) has to be finite for $\varepsilon:=\delta/2$, where
$\delta
:= d(f, \mathcal{\tilde{A}})$ is obviously positive since $\mathcal
{\tilde{A}}$ is closed.

With this choice of $\varepsilon$, we clearly have
%
%
\begin{equation}
\label{ex11} \PP\bigl\{\bigl\|\Gamma'_n W_{(n)}/c_n
- f\bigr\| < \varepsilon\bigr\} \le\PP\bigl\{d\bigl(\Gamma'_n
W_{(n)}/c_n, \mathcal{\tilde{A}}\bigr) \ge\varepsilon\bigr
\}.
\end{equation}
Define for $k \ge1$,
\begin{eqnarray*}
J_{k,\ell} &=& \bigl\{n\dvtx \exp(m_{k,\ell}) \le c_n
< \exp(n_{k,\ell})\bigr\},\qquad 0 \le\ell\le k,
\\
J'_{k,\ell} &=& \bigl\{n\dvtx \exp(n_{k,\ell}) \le
c_n < \exp(m_{k,\ell
+1})\bigr\},\qquad 0 \le\ell\le k.
\end{eqnarray*}
Employing once more inequality (\ref{exa2}) and recalling the
definition of $H$, we see that for all $n \in J_{k,\ell}$,
%
%
\begin{equation}
\label{ex12} \bigl|\Gamma'_n z_{\ell}\bigr|^2
=\E\bigl[\langle X, z_{\ell}\rangle^2 I\bigl\{|Z| \le
c_n\bigr\} \bigr]\le\bigl( \sigma_{\ell}^2
+2^{-k}\bigr)\log m_{k,l}
\end{equation}
and for any vector $w$ such that $\langle z_{\ell}, w \rangle=0, |w|=1$,
%
%
\begin{equation}
\label{ex13} \bigl|\Gamma'_n w\bigr|^2 = \E\bigl[
\langle X, w \rangle^2 I\bigl\{|Z| \le c_n\bigr\}\bigr] \le
2^{-k} \log m_{k,\ell}.
\end{equation}
Let again $\{w_{\ell,1},\ldots, w_{\ell,d}\}$ be an orthonormal basis
of $\R^d$ with $w_{\ell,1}=z_{\ell}$. Then,
\[
\Gamma'_n W_{(n)} = \sum
_{i=1}^d \bigl\langle w_{\ell,i},
\Gamma'_n W_{(n)}\bigr\rangle
w_{\ell, i}= \bigl\langle\Gamma'_n
z_{\ell}, W_{(n)}\bigr\rangle z_{\ell} + \sum
_{i=2}^d \bigl\langle\Gamma'_n
w_{\ell,i}, W_{(n)}\bigr\rangle w_{\ell, i}.
\]
Set $\varepsilon_1 =\varepsilon/\sqrt{d}$. Using the above representation and
the trivial fact that $d(f,\mathcal{\tilde{A}}) \le d(f,\mathcal
{K}_{\ell})$ for any $\ell\ge1$ with $\mathcal{K}_{\ell}=\sigma
_{\ell
}z_{\ell}\mathcal{K}$, we can infer from (\ref{ex11}) that
\begin{eqnarray*}
\PP\bigl\{\bigl\|\Gamma'_n W_{(n)}/c_n
- f\bigr\| < \varepsilon\bigr\}&\le& \PP\bigl\{ d\bigl(\bigl\langle
\Gamma'_n z_{\ell}, W_{(n)}\bigr
\rangle/c_n, \sigma_{\ell}\mathcal{K}\bigr) \ge
\varepsilon_1\bigr\}
\\
&&{} +\sum_{i=2}^d \PP\bigl\{\bigl\|\bigl
\langle\Gamma'_n w_{\ell,i}, W_{(n)}
\bigr\rangle\bigr\| \ge\varepsilon_1 c_n\bigr\}.
\end{eqnarray*}
To bound the first probability on the right-hand side, we make use of
an inequality for standard 1-dimensional Brownian motion $W'(t), 0 \le
t \le1$ which is due to Talagrand \cite{tala} and which states that
there exists an absolute constant $C>0$ such that for any $\lambda> 0$
and $x>0$,
\[
\PP\bigl\{d\bigl(W', \lambda\mathcal{K}\bigr) \ge x\bigr\} \le\exp
\bigl(Cx^{-2} - x\lambda /2 - \lambda^2/2\bigr).
\]
As $\langle\Gamma'_n z_{\ell}, W_{(n)}\rangle/\sqrt{n} \stackrel
{d}{=}|\Gamma'_n z_{\ell}|W'$, we have
\[
\PP\bigl\{d\bigl(\bigl\langle\Gamma'_n
z_{\ell}, W_{(n)}\bigr\rangle/c_n, \sigma
_{\ell
}\mathcal{K}\bigr) \ge\varepsilon_1\bigr\} =\PP\bigl\{d
\bigl(W', \lambda_n \mathcal{K}\bigr) \ge
x_n\bigr\},
\]
where $\lambda_n = \sigma_{\ell}\sqrt{2}\log\log n/|\Gamma'_n
z_{\ell
}|$ and $x_n = \varepsilon_1 \sqrt{2}\log\log n/|\Gamma'_n z_{\ell
}|$.

Recalling (\ref{ex12}) and noticing that $\log m_{k,\ell} \le\log
\log
c_n \le\log\log n, n \in J_{k,\ell}$, we can conclude that for $n
\in
J_{k,\ell}$ and large enough $k$,
\begin{eqnarray*}
&&\PP\bigl\{d\bigl(\bigl\langle\Gamma'_n
z_{\ell}, W_{(n)}\bigr\rangle/c_n, \sigma
_{\ell
}\mathcal{K}\bigr) \ge\varepsilon_1\bigr\}
\\
&&\qquad\le\exp \biggl(C\varepsilon_1^{-2}\frac{\sigma_{\ell}^2 + 2^{-k}}{2
\log
\log n} -
\frac{\sigma_{\ell}(\sigma_{\ell} + \varepsilon_1)}{\sigma
_{\ell
}^2 + 2^{-k}}\log\log n \biggr)
\\
&&\qquad\le 2(\log n)^{-1- \varepsilon_1/2}.
\end{eqnarray*}
To see this observe that if $n \in J_{k,\ell}$ and $1 \ge\sigma
_{\ell
}^2 \ge1/\ell, \ell\ge1$, we have for $\ell\le k+1$ and large $k$ that
\[
\frac{\sigma_{\ell}(\sigma_{\ell} + \varepsilon_1)}{\sigma_{\ell
}^2 +
2^{-k}} \ge\frac{1+ \varepsilon_1}{1 + 2^{-k}\sigma_{\ell}^{-2}} \ge \frac
{1+ \varepsilon_1}{1 + 2^{-k}(k+1)} \ge1 +
\varepsilon_1/2.
\]
Similarly, using the fact that $\langle\Gamma'_n w_{\ell,i},
W_{(n)}\rangle/\sqrt{n} \stackrel{d}{=} |\Gamma'_n w_{\ell,i}| W'$ in
conjunction with the inequality $\PP\{\|W'\| \ge x\} \le2\exp(-x^2/2),
x>0$, we get from (\ref{ex13}) for $n \in J_{k,\ell}$ and $2 \le i
\le
d$ that
\begin{eqnarray*}
&&\PP\bigl\{\bigl\|\bigl\langle\Gamma'_n w_{\ell,i},
W_{(n)}\bigr\rangle\bigr\| \ge \varepsilon_1 c_n\bigr
\}
\\
&&\qquad\le \PP\bigl\{\bigl|\Gamma'_n w_{\ell,i}\bigr|
\bigl\|W'\bigr\| \ge\varepsilon_1 \sqrt{2} \log \log n\bigr\}
\le2(\log n)^{-2^k \varepsilon_1^2}.
\end{eqnarray*}
It follows that
%
%
\begin{equation}
\label{ex14} \sum_{k=1}^{\infty} \sum
_{\ell=1}^k \sum_{n \in J_{k,\ell}}
n^{-1} \PP\bigl\{ \bigl\|\Gamma'_n
W_{(n)}/c_n - f\bigr\| < \varepsilon\bigr\} < \infty.
\end{equation}
We still need that
%
%
\begin{equation}
\label{ex15} \sum_{k=1}^{\infty} \sum
_{\ell=1}^k \sum_{n \in J'_{k,\ell}}
n^{-1} \PP \bigl\{\bigl\|\Gamma'_n
W_{(n)}/c_n - f\bigr\| < \varepsilon\bigr\} < \infty.
\end{equation}
To prove that, we simply note that $\|f\| \ge2\varepsilon$ since $d(f,
\tilde{\mathcal{A}}) = 2\varepsilon$ and
$0 \in\tilde{\mathcal{A}}$. Consequently, we have for any $n \ge1$,
\[
\PP\bigl\{\bigl\|\Gamma'_n W_{(n)}/c_n
- f\bigr\| < \varepsilon\bigr\} \le\PP\bigl\{\bigl\|\Gamma'_n
W_{(n)}\bigr\| \ge\varepsilon c_n\bigr\} \le\PP\bigl\{\bigl\|
\Gamma'_n\bigr\|\| W_{(n)}\| \ge\varepsilon
c_n\bigr\}.
\]
From the definition of $X$ it immediately follows that $\|\Gamma'_n\|^2
\le H(c_n) \le\log\log n,  n \ge3$.
Thus we have
\[
\PP\bigl\{\bigl\|\Gamma'_n W_{(n)}/c_n
- f\bigr\| < \varepsilon\bigr\} \le\PP\bigl\{\| W_{(n)}\| \ge\varepsilon\sqrt{2n
\log\log n}\bigr\}\le2d (\log n)^{-\varepsilon^2/d}.
\]
Using a similar argument as in the proof of (\ref{ex2}) and
$c_n^2/c_m^2 \ge n/m, m \le n$, we find that
\[
\sum_{n \in J'_{k,\ell}} n^{-1} \le2(m_{k,\ell+1}
- n_{k,\ell}) +1 = 2k^3 +1,\qquad 1 \le\ell\le k, k \ge1.
\]
As we have $\log m_{k,0} \ge2^{k^3}, k \ge1$, we can conclude that
\[
\sum_{\ell=1}^k \sum
_{n \in J'_{k,\ell}} n^{-1} \PP\bigl\{\bigl\|\Gamma'_n
W_{(n)}/c_n - f\bigr\| < \varepsilon\bigr\} \le2d
\bigl(2k^4 + k\bigr)2^{-\varepsilon^2 k^3/d},\qquad k \ge1,
\]
which trivially implies (\ref{ex15}).

Combining (\ref{ex14}) and (\ref{ex15}), it follows from (\ref{ex10})
that $f \notin C(\{s_n\})$. We see that (\ref{ex3}) holds and the
theorem has been proven.
\end{pf*}
%

%



\printaddresses


\begin{thebibliography}{20}

\bibitem{Day}
%
\begin{bbook}[mr]
\bauthor{\bsnm{Day},~\bfnm{Mahlon~M.}\binits{M.~M.}}
(\byear{1973}).
\btitle{Normed Linear Spaces},
\bedition{3rd} ed.
\bpublisher{Springer}, \blocation{New York}.
\bid{mr={0344849}}%
\bptok{imsref}%
\end{bbook}
%
\endbibitem

\bibitem{dlpG}
%
\begin{bbook}[mr]
\bauthor{\bparticle{de~la}
\bsnm{Pe{\~n}a},~\bfnm{V{\'{\i}}ctor~H.}\binits{V.~H.}} \AND
\bauthor{\bsnm{Gin{\'e}},~\bfnm{Evarist}\binits{E.}}
(\byear{1999}).
\btitle{Decoupling: From Dependence to Independence, Randomly Stopped
Processes.
$U$-Statistics and Processes. Martingales and Beyond}.
\bpublisher{Springer}, \blocation{New York}.
\bid{doi={10.1007/978-1-4612-0537-1}, mr={1666908}}
\bptok{imsref}%
\end{bbook}
%
\endbibitem

\bibitem{E-3}
%
\begin{barticle}[mr]
\bauthor{\bsnm{Einmahl},~\bfnm{Uwe}\binits{U.}}
(\byear{1993}).
\btitle{Toward a general law of the iterated logarithm in {B}anach space}.
\bjournal{Ann. Probab.}
\bvolume{21}
\bpages{2012--2045}.
\bid{issn={0091-1798}, mr={1245299}}
\bptok{imsref}%
\end{barticle}
%
\endbibitem

\bibitem{E-1}
%
\begin{barticle}[mr]
\bauthor{\bsnm{Einmahl},~\bfnm{Uwe}\binits{U.}}
(\byear{1995}).
\btitle{On the cluster set problem for the generalized law of the iterated
logarithm in {E}uclidean space}.
\bjournal{Ann. Probab.}
\bvolume{23}
\bpages{817--851}.
\bid{issn={0091-1798}, mr={1334174}}
\bptok{imsref}%
\end{barticle}
%
\endbibitem

\bibitem{E}
%
\begin{barticle}[mr]
\bauthor{\bsnm{Einmahl},~\bfnm{Uwe}\binits{U.}}
(\byear{2007}).
\btitle{A generalization of {S}trassen's functional {LIL}}.
\bjournal{J. Theoret. Probab.}
\bvolume{20}
\bpages{901--915}.
\bid{doi={10.1007/s10959-007-0091-0}, issn={0894-9840}, mr={2359061}}
\bptok{imsref}%
\end{barticle}
%
\endbibitem

\bibitem{E-2}
%
\begin{barticle}[mr]
\bauthor{\bsnm{Einmahl},~\bfnm{U.}\binits{U.}}
(\byear{2009}).
\btitle{A new strong invariance principle for sums of independent random
vectors}.
\bjournal{J. Math. Sci.}
\bvolume{163}
\bpages{311--327}.
\bptok{imsref}%
\end{barticle}
%
\endbibitem

\bibitem{EK}
%
\begin{barticle}[mr]
\bauthor{\bsnm{Einmahl},~\bfnm{U.}\binits{U.}} \AND
\bauthor{\bsnm{Kuelbs},~\bfnm{J.}\binits{J.}}
(\byear{2001}).
\btitle{Cluster sets for a generalized law of the iterated logarithm in
{B}anach spaces}.
\bjournal{Ann. Probab.}
\bvolume{29}
\bpages{1451--1475}.
\bid{doi={10.1214/aop/1015345758}, issn={0091-1798}, mr={1880228}}
\bptok{imsref}%
\end{barticle}
%
\endbibitem

\bibitem{EL-1}
%
\begin{barticle}[mr]
\bauthor{\bsnm{Einmahl},~\bfnm{Uwe}\binits{U.}} \AND
\bauthor{\bsnm{Li},~\bfnm{Deli}\binits{D.}}
(\byear{2005}).
\btitle{Some results on two-sided {LIL} behavior}.
\bjournal{Ann. Probab.}
\bvolume{33}
\bpages{1601--1624}.
\bid{doi={10.1214/009117905000000198}, issn={0091-1798}, mr={2150200}}
\bptok{imsref}%
\end{barticle}
%
\endbibitem

\bibitem{EL}
%
\begin{barticle}[mr]
\bauthor{\bsnm{Einmahl},~\bfnm{Uwe}\binits{U.}} \AND
\bauthor{\bsnm{Li},~\bfnm{Deli}\binits{D.}}
(\byear{2008}).
\btitle{Characterization of {LIL} behavior in {B}anach space}.
\bjournal{Trans. Amer. Math. Soc.}
\bvolume{360}
\bpages{6677--6693}.
\bid{doi={10.1090/S0002-9947-08-04522-4}, issn={0002-9947}, mr={2434306}}
\bptok{imsref}%
\end{barticle}
%
\endbibitem

\bibitem{gr}
%
\begin{barticle}[mr]
\bauthor{\bsnm{Grill},~\bfnm{Karl}\binits{K.}}
(\byear{1991}).
\btitle{A {$\liminf$} result in {S}trassen's law of the iterated logarithm}.
\bjournal{Probab. Theory Related Fields}
\bvolume{89}
\bpages{149--157}.
\bid{doi={10.1007/BF01366903}, issn={0178-8051}, mr={1110535}}
\bptok{imsref}%
\end{barticle}
%
\endbibitem

\bibitem{Ke}
%
\begin{barticle}[mr]
\bauthor{\bsnm{Kesten},~\bfnm{Harry}\binits{H.}}
(\byear{1970}).
\btitle{The limit points of a normalized random walk}.
\bjournal{Ann. Math. Statist.}
\bvolume{41}
\bpages{1173--1205}.
\bid{issn={0003-4851}, mr={0266315}}
\bptok{imsref}%
\end{barticle}
%
\endbibitem

\bibitem{Kla}
%
\begin{barticle}[mr]
\bauthor{\bsnm{Klass},~\bfnm{M.~J.}\binits{M.~J.}}
(\byear{1976}).
\btitle{Toward a universal law of the iterated logarithm. {I}}.
\bjournal{Z. Wahrsch. Verw. Gebiete}
\bvolume{36}
\bpages{165--178}.
\bid{mr={0415742}}
\bptok{imsref}%
\end{barticle}
%
\endbibitem

\bibitem{Kue}
%
\begin{barticle}[mr]
\bauthor{\bsnm{Kuelbs},~\bfnm{J.}\binits{J.}}
(\byear{1981}).
\btitle{When is the cluster set of {$S_{n}/a_{n}$} empty?}
\bjournal{Ann. Probab.}
\bvolume{9}
\bpages{377--394}.
\bid{issn={0091-1798}, mr={0614624}}
\bptok{imsref}%
\end{barticle}
%
\endbibitem

\bibitem{Kue2}
%
\begin{barticle}[mr]
\bauthor{\bsnm{Kuelbs},~\bfnm{J.}\binits{J.}}
(\byear{1985}).
\btitle{The {LIL} when {$X$} is in the domain of attraction of a {G}aussian
law}.
\bjournal{Ann. Probab.}
\bvolume{13}
\bpages{825--859}.
\bid{issn={0091-1798}, mr={0799424}}
\bptok{imsref}%
\end{barticle}
%
\endbibitem

\bibitem{KLL}
%
\begin{barticle}[mr]
\bauthor{\bsnm{Kuelbs},~\bfnm{James}\binits{J.}},
\bauthor{\bsnm{Li},~\bfnm{Wenbo~V.}\binits{W.~V.}} \AND
\bauthor{\bsnm{Linde},~\bfnm{Werner}\binits{W.}}
(\byear{1994}).
\btitle{The {G}aussian measure of shifted balls}.
\bjournal{Probab. Theory Related Fields}
\bvolume{98}
\bpages{143--162}.
\bid{doi={10.1007/BF01192511}, issn={0178-8051}, mr={1258983}}
\bptok{imsref}%
\end{barticle}
%
\endbibitem

\bibitem{kzi}
%
\begin{barticle}[mr]
\bauthor{\bsnm{Kuelbs},~\bfnm{Jim}\binits{J.}} \AND
\bauthor{\bsnm{Zinn},~\bfnm{Joel}\binits{J.}}
(\byear{2008}).
\btitle{Another view of the {CLT} in {B}anach spaces}.
\bjournal{J.~Theoret. Probab.}
\bvolume{21}
\bpages{982--1029}.
\bid{doi={10.1007/s10959-008-0166-6}, issn={0894-9840}, mr={2443644}}
\bptok{imsref}%
\end{barticle}
%
\endbibitem

\bibitem{LT}
%
\begin{bbook}[mr]
\bauthor{\bsnm{Lindenstrauss},~\bfnm{Joram}\binits{J.}} \AND
\bauthor{\bsnm{Tzafriri},~\bfnm{Lior}\binits{L.}}
(\byear{1977}).
\btitle{Classical {B}anach Spaces. {I}: Sequence Spaces}.
\bpublisher{Springer}, \blocation{Berlin}.
\bid{mr={0500056}}
\bptok{imsref}%
\end{bbook}
%
\endbibitem

\bibitem{sakh}
%
\begin{bincollection}[mr]
\bauthor{\bsnm{Sakhanenko},~\bfnm{Alexander~I.}\binits{A.~I.}}
(\byear{2000}).
\btitle{A new way to obtain estimates in the invariance principle}.
In \bbooktitle{High Dimensional Probability, {II} ({S}eattle, {WA}, 1999)}.
\bseries{Progress in Probability}
\bvolume{47}
\bpages{223--245}.
\bpublisher{Birkh\"auser}, \blocation{Boston, MA}.
\bid{mr={1857325}}
\bptok{imsref}%
\end{bincollection}
%
\endbibitem

\bibitem{tala}
%
\begin{bincollection}[auto:STB|2013/09/19|12:14:10]
\bauthor{\bsnm{Talagrand},~\bfnm{M.}\binits{M.}}
(\byear{1992}).
\btitle{On the rate of clustering in Strassen's LIL for Brownian motion}.
In \bbooktitle{Probability in Banach Spaces}
(\beditor{\bfnm{R.}\binits{R.}~\bsnm{Dudley}},
\beditor{\bfnm{M.}\binits{M.}~\bsnm{Hahn}} \AND
\beditor{\bfnm{J.}\binits{J.}~\bsnm{Kuelbs}}, eds.)
\bvolume{8}
\bpages{333--347}.
\bpublisher{Birkh\"auser}, \blocation{Boston}.
\bptok{imsref}%
\end{bincollection}
%
\endbibitem

\end{thebibliography}
\end{document}